\documentclass[onecolumn,final,10pt]{elsarticle}
\usepackage{hyperref}
\usepackage{url}
\usepackage{caption}
\usepackage{subcaption}
\usepackage{units}
\usepackage{mathematix}
\usepackage{listings}
\usepackage{lipsum}
\usepackage{eurosym}
\usepackage{microtype}
\usepackage{geometry}

\usepackage{flushend}

\colorlet{punct}{red!60!black}
\definecolor{background}{HTML}{EEEEEE}
\definecolor{delim}{RGB}{20,105,176}
\colorlet{numb}{magenta!60!black}

\lstdefinelanguage{json}{
    basicstyle=\normalfont\ttfamily,
    stepnumber=1,
    numbersep=8pt,
    showstringspaces=false,
    breaklines=true,
    frame=lines,
    backgroundcolor=\color{background},
    literate=
     *{0}{{{\color{numb}0}}}{1}
      {1}{{{\color{numb}1}}}{1}
      {2}{{{\color{numb}2}}}{1}
      {3}{{{\color{numb}3}}}{1}
      {4}{{{\color{numb}4}}}{1}
      {5}{{{\color{numb}5}}}{1}
      {6}{{{\color{numb}6}}}{1}
      {7}{{{\color{numb}7}}}{1}
      {8}{{{\color{numb}8}}}{1}
      {9}{{{\color{numb}9}}}{1}
      {:}{{{\color{punct}{:}}}}{1}
      {,}{{{\color{punct}{,}}}}{1}
      {\{}{{{\color{delim}{\{}}}}{1}
      {\}}{{{\color{delim}{\}}}}}{1}
      {[}{{{\color{delim}{[}}}}{1}
      {]}{{{\color{delim}{]}}}}{1},
}

\hypersetup{
    unicode=false,          % non-Latin characters in Acrobat’s bookmarks
    pdftoolbar=true,        % show Acrobat’s toolbar?
    pdfmenubar=true,        % show Acrobat’s menu?
    pdffitwindow=false,     % window fit to page when opened
    pdfstartview={FitH},    % fits the width of the page to the window
    pdfnewwindow=true,      % links in new PDF window
    colorlinks=true,       % false: boxed links; true: colored links
    linkcolor=red,          % color of internal links (change box color with linkbordercolor)
    citecolor=green,        % color of links to bibliography
    filecolor=magenta,      % color of file links
    urlcolor=cyan           % color of external links
}

\newcommand{\dist}{\operatorname{dist}}
\journal{Environmental Modelling \& Software}
\def\name{RapidNet}

\begin{document}

\begin{frontmatter}

%% Title, authors and addresses

%% use the tnoteref command within \title for footnotes;
%% use the tnotetext command for theassociated footnote;
%% use the fnref command within \author or \address for footnotes;
%% use the fntext command for theassociated footnote;
%% use the corref command within \author for corresponding author footnotes;
%% use the cortext command for theassociated footnote;
%% use the ead command for the email address,
%% and the form \ead[url] for the home page:
%% \title{Title\tnoteref{label1}}
%% \tnotetext[label1]{}
%% \author{Name\corref{cor1}\fnref{label2}}
%% \ead{email address}
%% \ead[url]{home page}
%% \fntext[label2]{}
%% \cortext[cor1]{}
%% \address{Address\fnref{label3}}
%% \fntext[label3]{}

\title{Uncertainty-aware demand management of water distribution networks in deregulated energy markets}

%% use optional labels to link authors explicitly to addresses:
%% \author[label1,label2]{}
%% \address[label1]{}
%% \address[label2]{}

\author[KUL]{P. Sopasakis\corref{cor1}}\ead{pantelis.sopasakis@kuleuven.be}
\author[TUB]{A.K. Sampathirao\fnref{fnajay}}\ead{sampathirao@control.tu-berlin.de}
\author[IMT]{A. Bemporad}\ead{alberto.bemporad@imtlucca.it}
\author[KUL]{P. Patrinos\fnref{fn1}}\ead{panos.patrinos@esat.kuleuven.be}

\cortext[cor1]{Corresponding author.}
\fntext[fnajay]{The work of A. K. Sampathirao was supported by the German Federal Ministry for Economic Affairs and 
                Energy (BMWi), Project No. 0324024A.}
\fntext[fn1]{The work of P. Patrinos was supported by the KU Leuven Research Council under BOF/STG-15-043.}
\address[KUL]{KU Leuven, Department of Electrical Engineering (ESAT),
              STADIUS Center for Dynamical Systems, Signal Processing and Data Analytics
              \& Optimization in Engineering (OPTEC),
              Kasteelpark Arenberg 10, 3001 Leuven, Belgium.}
\address[TUB]{Technische Universit\"at Berlin, Fachgebiet Regelungssysteme, Einsteinufer 17, D-10587 Berlin, Germany}
\address[IMT]{IMT School for Advanced Studies Lucca,
                   Piazza San Ponziano 6,
                   Lucca 55100, Italy.}

\begin{abstract}
%% Text of abstract
We present an open-source solution for the operational control of drinking 
water distribution networks which accounts for the inherent uncertainty 
in water demand and electricity prices in the day-ahead market of a volatile 
deregulated economy. As increasingly more energy markets adopt this trading 
scheme, the operation of drinking water networks requires uncertainty-aware 
control approaches that mitigate the effect of volatility and result in an 
economic and safe operation of the network that meets the consumers' need 
for uninterrupted water supply. We propose the use of scenario-based stochastic
model predictive control: an advanced control methodology which comes at a 
considerable computation cost which is overcome by harnessing the parallelization
capabilities of graphics processing units (GPUs) and using a massively 
parallelizable algorithm based on the accelerated proximal gradient method.
\end{abstract}

\begin{keyword}
Drinking water networks \sep
Bid-based energy market \sep
Stochastic model predictive control \sep 
Graphics processing units \sep 
Scenario trees \sep 
Open-source software
%% keywords here, in the form: keyword \sep keyword
%
%% PACS codes here, in the form: \PACS code \sep code
%
%% MSC codes here, in the form: \MSC code \sep code
%% or \MSC[2008] code \sep code (2000 is the default)
\end{keyword}

\end{frontmatter}

%% \linenumbers
%% main text

\section*{Software availability}

\noindent %
\emph{Name of software:} \texttt{\name}.\\[0.35em]
\emph{Hardware required:} CUDA-compliant GPU (tested on NVIDIA Tesla 2075).\\[0.35em]
\emph{Software required:} CUDA framework v6.0 or higher (including cuBLAS) and rapidjson (\url{http://rapidjson.org/}).\\[0.35em]
\emph{Availability:} Open-source software, licence: GNU LGPL v3.0. Available online at \url{https://github.com/GPUEngineering/RapidNet}.\\[0.35em]
\emph{Program language}: \texttt{CUDA-C++}/\texttt{C++}.\\[0.35em]
\emph{First release:} 2017.

\section{Introduction}\label{sec:into}

\subsection{State of the art}
% Smart systems - about (cyber-physical systems)
The explosive proliferation of interconnected sensing, computing and
communication devices has marked the advent of the concept of
\emph{cyber-physical systems} --- ensembles of computational
and physical components. In drinking water networks this trend
has ushered in new control paradigms where the profusion of
data, produced by a network of sensors and stored in a database,
is used to prescribe informed control actions%
~\citep{EggMutWan+17,Meseguer2017,Solomatine03,LobSol2002}.
Nevertheless, as these data, be they water demand values or 
electricity prices, cannot be modeled perfectly, the associated
uncertainty is shifted to the decision making process.

% Introduction and problem statement
% --> Opening
The high uncertainty in the operation of drinking water networks,
as a result of the volatility of future demands as well as energy
prices (in a deregulated energy market) is likely to lead to a rather
expensive operating mode with poor quality of service (the network
may not always be able to provide the necessary amount of water to the
consumers).
In control engineering practice, this uncertainty is often addressed in
a worst-case fashion~\citep{SamGroSop+14,WanOcaPui15} --- if not neglected at all
--- leading to conservative and suboptimal control policies.
% --> Objective...
It is evident that it is necessary to devise control methods which take into
account the probabilistic nature of the underlying uncertainty making use of the
wealth of available historical data aiming at a proactive and foresightful control
scheme which leads to an improved closed-loop performance. These requirements
necessitate the use of \textit{stochastic model predictive control}: an
advanced control methodology where at every time instant we determine a
sequence of \textit{control laws} which minimizes the \textit{expected value}
of a performance index taken with respect to the distribution of the
uncertainty~\citep{Mesbah2016}.
Optimization-based approaches for the operational management of 
water networks have been studied and are well established in engineering practice~\cite{MalaJetmarova2017209}.

% Our methodology
Indeed, scenario-based stochastic model predictive control (SSMPC) has been shown
to lead to remarkable decrease in the operating cost and improvement in the
quality of service of drinking water networks~\citep{SamSopBemPat17b}.
In SSMPC, the uncertain disturbances are treated as random variables on
a discrete sample space without assuming any parametric form for their 
distribution~\citep{CalCam06}.
The scenario approach was identified in a recent review as a powerful
method for mitigating uncertainty in environmental modeling related to 
water management~\cite{Horne2016326}.
Although this approach offers a realistic control solution as it is 
entirely data-driven, this comes with considerable computational burden as the 
resulting optimization problems are of particularly large 
scale~\citep{GorNEm14,GrossoTrees14}. This has rendered the use of
SSMPC prohibitive and has hindered its applicability.
Indeed, hitherto there have been used only conventional model predictive
control approaches~\citep{OcaFamBarPui10,BaVrePa+13}, robust worst-case formulations%
~\citep{SamGroSop+14,OcPuCe+09,TraBrd09}
and stochastic formulations where the underlying uncertainty is assumed
to be normally identically independently distributed~\citep{WanOcaPui16,GroVelOca+16}.
Note that it has been observed that demand prediction errors are 
typically follow heavy-tail distributions which cannot be well approximated 
by normal ones~\cite{Hutton201587}.

In this paper, we present a software for the fast and efficient solution
of such problems harnessing the immense computational capabilities of
graphics processing units (GPU) building up on our previous work%
~\citep{SamSopBemPat17b,SamSopBemPat16,SamSopBemPat15}.

% The SoTA in algorithms (how do people solve stochastic optimal control
% problems nowadays)?
There has been recently a lot of interest in the development of 
efficient methods for stochastic optimal control problems such as 
stochastic gradient methods~\citep{TheVilPatBem16}, the alternating 
directions method of multipliers (ADMM)~\citep{KanRagDiCai16} and
various decomposition methods which can lead to parallelizable 
methods~\citep{carpentier:hal-01232179,Defourny12} (the most popular being
the stochastic dual approximate dynamic programming~\citep{JiaPhaPow+14},
the progressive hedging approach~\citep{CarGenBas13} and dynamic
programming~\citep{Bert00}).
There have been proposed parallelizable interior point algorithms for
two-stage stochastic optimal control problems such as%
~\citep{KlinGros17,lubin-scalable-2011,jia2014,ChiPetZav14}
and an \textit{ad hoc} interior point solver for multi-stage problems~\citep{hubner2016distributed}.
However, interior point algorithms involve complex steps and are not suitable
for an implementation on GPUs which can make the most of the capabilities 
of the hardware.
Additionally, interior point methods cannot accommodate complex non-quadratic 
terms in the cost function such as soft constraints (distance-to-set functions).

% Why is this software necessary - competitors?
At large, not many software and libraries are available for stochastic
optimal control; one of the very few one may find on the web is
\href{http://gwr3n.github.io/jsdp/}{JSPD}, a generic Java stochastic dynamic
programming library.
\href{http://www.quantego.com/}{QUASAR} is a commercial tool for
scenario-based stochastic optimization.
One of the most popular tools in the toolbox of the water networks
engineer is PLIO~\citep{plio}, which implements MPC
algorithms. This work covers the yawning gap between engineering 
practice and the latest  developments in control and optimization 
theory for drinking water networks. 
These results can also be applied for the control
of other infrastructure with similar structure such as 
power grids~\citep{HanSopBemRaiCol15}.

\subsection{Contributions and Novelty}
% The algorithm
Despite the fact that SSMPC problems typically involve millions of decision variables,
the associated optimization problems possess a rich structure which can be exploited
to devise parallelizable \emph{ad hoc} methods to solve the problem more than an
order of magnitude faster than commercial solvers running on CPU.

% Software architecture
The architecture of our implementation comprises three independent
modules: (i) the \emph{network} module, (ii) the \emph{energy prices
and water demands forecaster} and (iii) the \emph{control module}.
The network module provides a dynamical system model which describes the
flow of water across the network together with the storage limits of the tanks
and the constraints on pumping capacities.
The network module defines a \emph{safety storage level} for each tank ---
a level which ensures the availability of water in case of high demand and
the maintenance of a minimum required pressure.
The forecasters produce a \emph{scenario tree}, that is, a tree of likely
future water demands and energy prices, upon which a contingency plan is made
by minimizing a cost function which quantifies the operating cost and the
quality of service. Such scenario trees are constructed from historical data
of energy prices and water demands.
The control module computes flow set-points, which are sent to the
pumping stations and valves, by solving a scenario-based stochastic model
predictive control problem over a finite prediction horizon.

% Assessment of closed-loop performance
The proposed stochastic model predictive controller leads to measurable benefits
for the operation of the water network. It leads to a more economic operation compared
to methods which do not take into consideration the stochastic nature of the energy
prices and water demands. In this paper, we assess the performance of the controlled
network using three key performance indicators: (i) the \emph{economic index},
(ii) the \emph{safety index}, which quantifies the extent of violation of the safety 
storage level requirement and (iii) the \emph{computational complexity index} with 
which we assess the computational feasibility of the controller.
Simulation results are provided using data from the water network of Barcelona and
the energy market of Austria.
The advantages of the adopted control methodology are combined with the computational 
power of GPUs, which enables us to solve problems of very large scale.

\subsection{Software}
% Software
Our implementation is available as an open-source and free software which can be readily
tailored to the needs of different water networks modifying the parameters of its modules.
The implementation is done entirely in \texttt{CUDA-C++} and it can be configured either
programmatically or using configuration files. The adopted object-oriented programming model
is amenable to extensions and users may specify their own predictive models, scenario trees,
cost functions, dynamical models and constraints. Our results are accompanied by extensive
benchmarks and the software is verified with unit tests.

\section{Modeling}

\subsection{Hydraulic modeling}
The water network involves four types of elements: water tanks,
active elements (pumps and valves), mixing nodes and demand sectors.
We focus on flow-based water networks where all flows are manipulated
variables based on the modeling approaches presented in%
~\citep{Batchabani2014,OcaFamBarPui10,Papageorgiou1984,GrossoTrees14,SamGroSop+14}.
\begin{figure*}
 \centering
 \newlength{\subfiglength}
 \setlength{\subfiglength}{0.19\linewidth}
 \begin{subfigure}[t]{\subfiglength}
		\centering
		\includegraphics[width=1.0\linewidth]{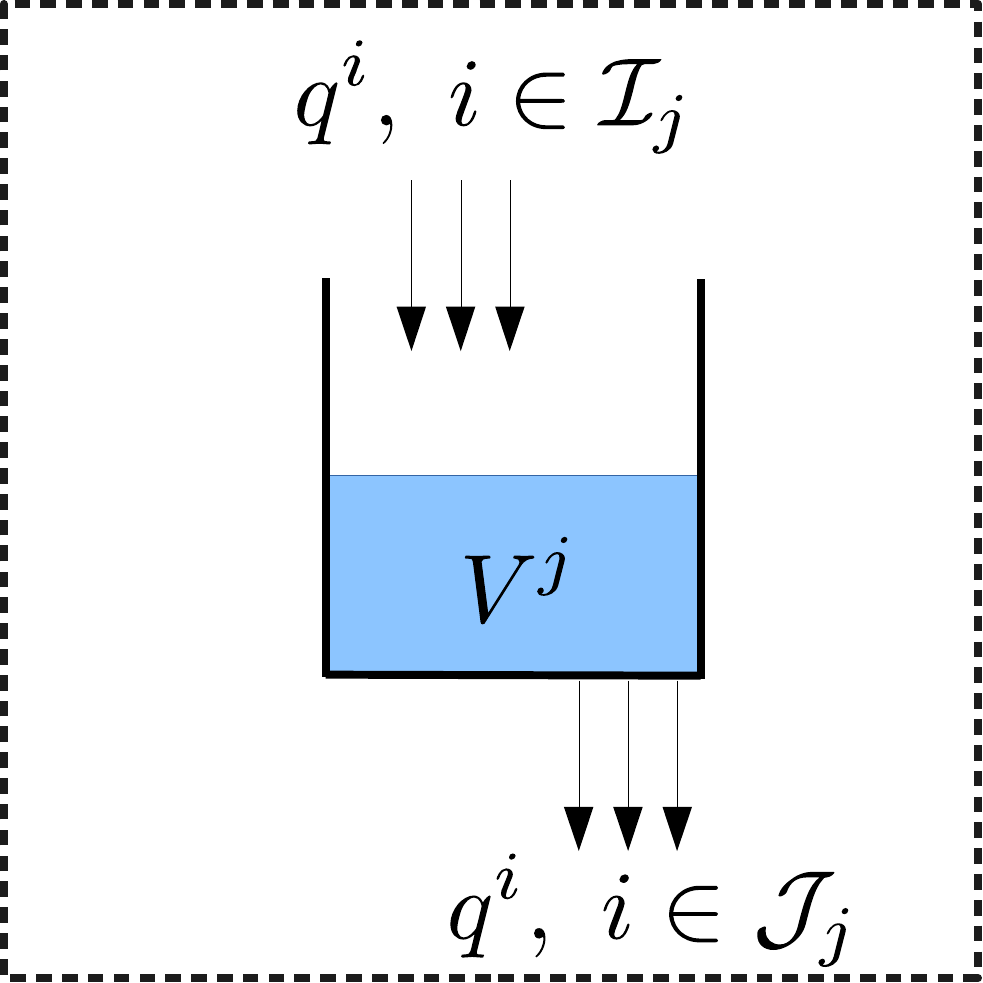}
		\caption{Water tank}\label{fig:dwn:a}
 \end{subfigure}
 \begin{subfigure}[t]{\subfiglength}
		\centering
		\includegraphics[width=1.0\linewidth]{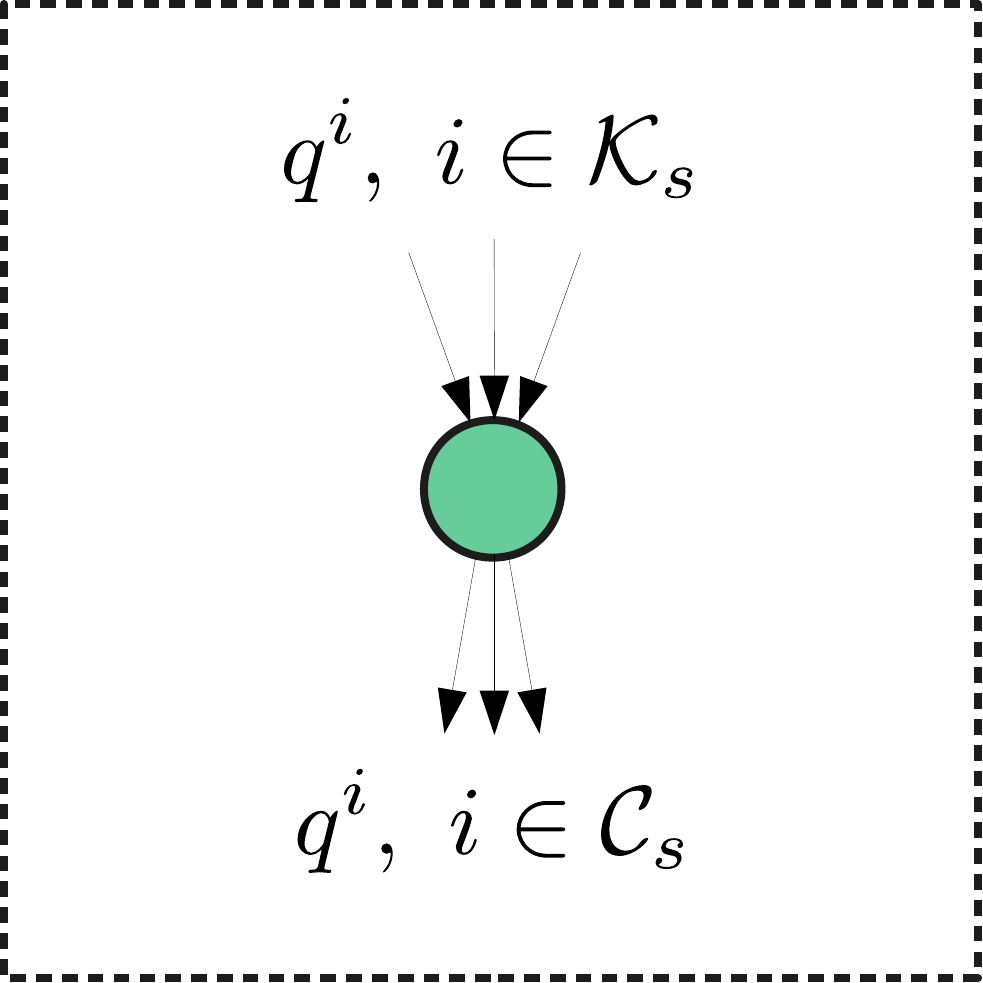}
		\caption{Mixing node}\label{fig:dwn:b}
 \end{subfigure}
 \begin{subfigure}[t]{\subfiglength}
		\centering
		\includegraphics[width=1.0\linewidth]{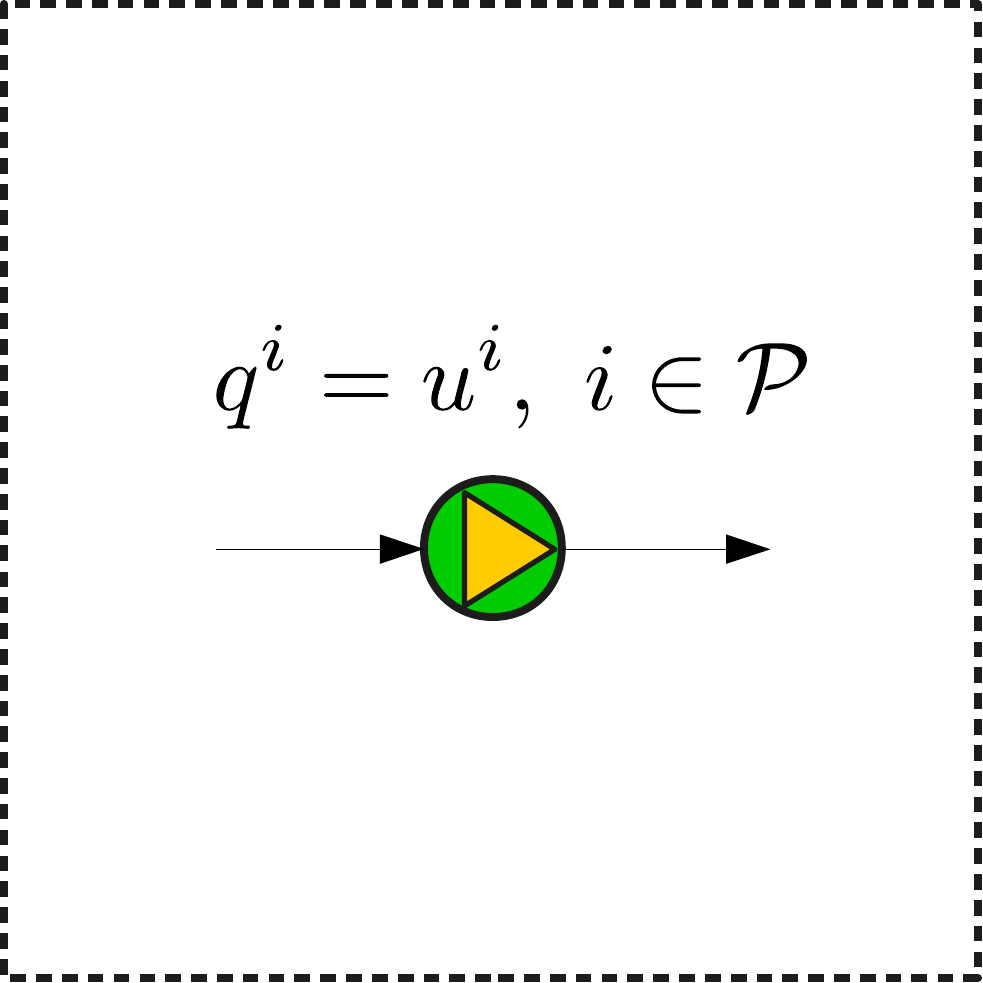}
		\caption{Controlled pump}\label{fig:dwn:c}
 \end{subfigure}
 \begin{subfigure}[t]{\subfiglength}
		\centering
		\includegraphics[width=1.0\linewidth]{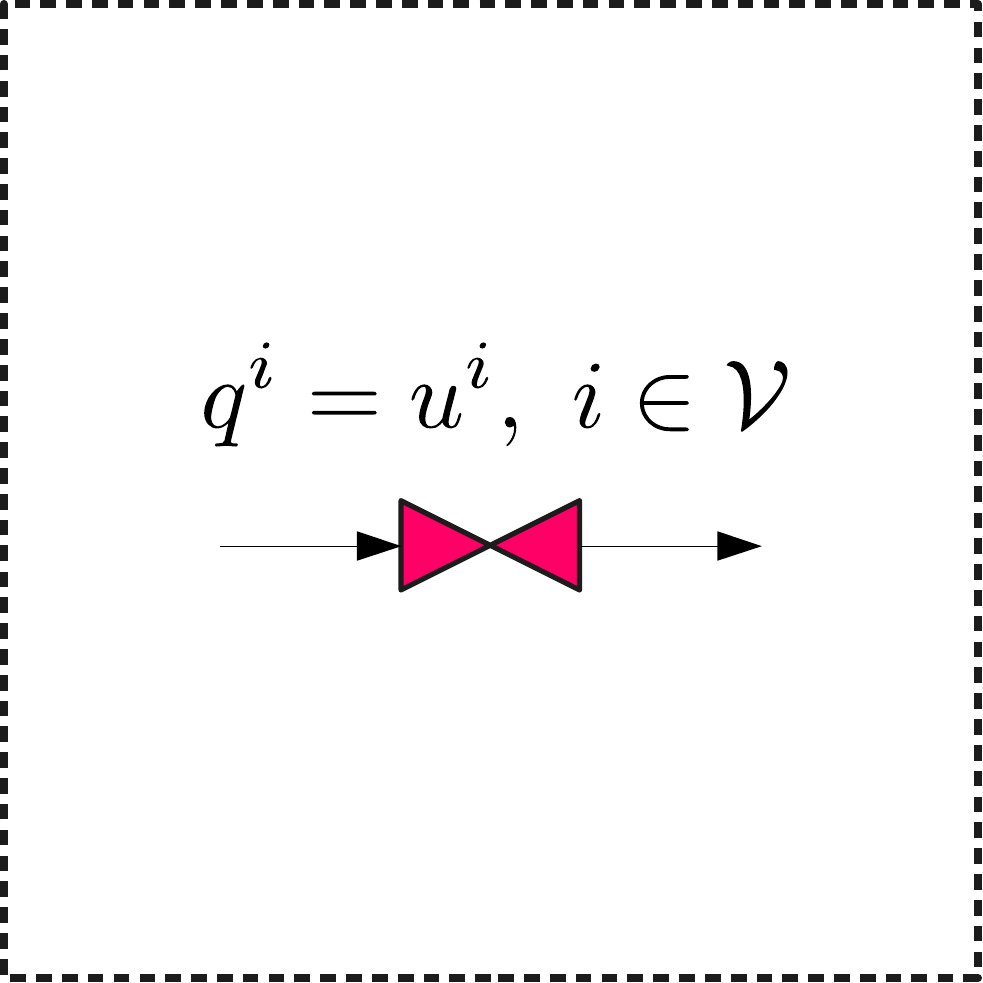}
		\caption{Controlled valve}\label{fig:dwn:d}
 \end{subfigure}
 \begin{subfigure}[t]{\subfiglength}
		\centering
		\includegraphics[width=1.0\linewidth]{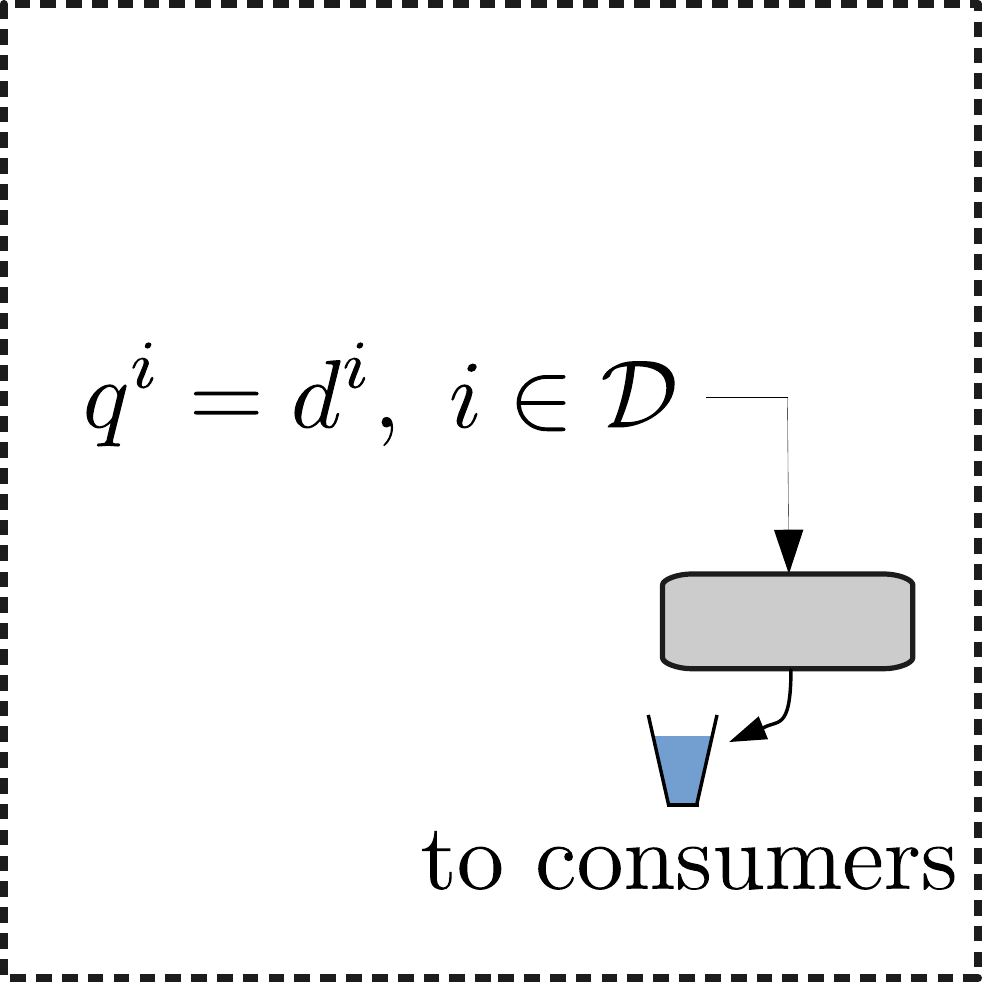}
		\caption{Demand sector}\label{fig:dwn:e}
 \end{subfigure}
 \caption{Elements of the water network and the associated flows.}
 \hrule
\end{figure*}

\textit{Water tanks} which play a crucial role in demand management as they ensure water supply,
obviate the need for continuous pumping and allow water to be pumped when the
price of electricity is lower and provide water in cases of unexpected peaks in
demand. Tank dynamics are modelled by simple mass balance equations: let $V^j(t)$
be the volume of tank $j=1,\ldots,N_t$ at time $t$.
We denote all flows in the network by $q^i$, $i=1,\ldots, N_f$.
Let $q^{i}(t)$, $i\in\mathcal{I}_j$ be the controlled inflowing streams to tank
$j$ and $q^{i}(t)$, $i\in\mathcal{J}_j$ be the controlled outflowing streams.
Then, the mass balance equation becomes (see Fig.~\ref{fig:dwn:a})
\begin{align}
 \label{eq:tank-dynamics}
 \frac{\d V^j(t)}{\d t} = \sum_{i\in\mathcal{I}_j}q^{i}(t) - \sum_{i\in\mathcal{J}_j}q^{i}(t).
\end{align}
The volume in each tank should never exceed a maximum limit $V^j_{\max}$
and it should always be above a hard lower limit $V^j_{\min}$, that is,
\begin{align}
 \label{eq:level-constraints}
 V^j_{\min} \leq V^j(t) \leq V^j_{\max}.
\end{align}

In addition, for the sake of service
reliability (availability of water when demand rises unexpectedly) and safety,
it is required that the level of water remains above a certain level $V^j_{\mathrm{safe}}$.
This is allowed to be violated occasionally, when the demand happens to be
too high.

\textit{Mixing nodes} are intersections where flows of water are merged
or separated. The mass balance equations for mixing nodes give rise to
algebraic constraints of the form
\begin{align}
 \sum_{i\in\mathcal{K}_s}q^i(t) = \sum_{i\in\mathcal{C}_s}q^i(t),
\end{align}
where $\mathcal{K}_s$ and $\mathcal{C}_s$ are the sets of indices of
the incoming and outgoing flows at node $s=1,\ldots, N_s$ as shown in
Fig.~\ref{fig:dwn:b}.

\textit{Pumps} and \textit{valves} are used to control the flow of water
in the network and transfer it across tanks and to the demand sectors.
We treat these as controlled systems --- indeed, pumping stations
and valves are equipped with local controllers --- to which we prescribe
flow set-points. The local control systems operate at a sampling rate of
about $\unit[1]{Hz}$, while the operational management of the network
updates its decisions at a much slower rate (e.g., hourly). It is reasonable
to assume that the local control system equilibrates fast enough to
neglect its dynamics in the context of operational control.
That said, the flow determined by each pump $i\in\mathcal{P}$ is equal to
its prescribed set-point $u^i$. As shown in Fig.~\ref{fig:dwn:c}, that is
\begin{align}
 q^i(t) = u^i(t),\ i\in\mathcal{P}.
\end{align}
Similarly, as shown in Fig.~\ref{fig:dwn:d} the flow through each valve is
\begin{align}
 q^i(t) = u^i(t),\ i\in\mathcal{V}.
\end{align}

All flows in the network are unidirectional, so we require that $q^i\geq 0$
for all $i=1,\ldots, N_f$. Each pump $i\in\mathcal{P}$ has a maximum
\textit{pumping capacity} $q^i_{\max}$, that is we require that
\begin{align}
 \label{eq:flow-constraints}
 0 \leq q^i(t)\leq q^i_{\max},\ i\in\mathcal{P}.
\end{align}

\textit{Demand sectors} are the exit nodes of the water from the
network towards the consumers as shown in Fig.~\ref{fig:dwn:e}.
At each demand sector $i\in\mathcal{D}$, the mass
balance yields
\begin{align}
 q^i(t) = d^i(t),\ i\in\mathcal{D}.
\end{align}
We treat $d^i(t)$ as a random process and elaborate on that in 
Section~\ref{sec:demand-electricity}.

We may now describe the system dynamics in terms of the state variable $x(t) = (V^j(t))_{j=1,\ldots,N_t}$,
the input variable $u(t)=(u^j(t))_{j\in\mathcal{P}\cup\mathcal{V}}$ 
and the disturbance $d(t) = (d^i(t))_{i\in\mathcal{D}}$.
{ 
Let $n_u=|\mathcal{P}\cup\mathcal{V}|$ be the dimension of $u(t)$, that is,
the total number of pumps and valves. Let $n_d=|\mathcal{D}|$
be the total number of demand nodes.}
By discretizing the dynamical equation~\eqref{eq:tank-dynamics} 
using the exact discretization method~\cite[Sec.~{4.2.1}]{tsong1999book}
and taking into account the algebraic constraints stated above, we may
write the system dynamics in the following form of a discrete-time
linear time-invariant system
\begin{subequations}
\begin{align}
 x_{k+1} &= Ax_{k} + Bu_{k} + G_d d_{k},\label{eq:dynamics}\\
 0 &= Eu_{k} + E_d d_{k}\label{eq:input-disturbance-coupling},
\end{align}
\end{subequations}
where {
$A\in\Re^{N_t\times N_t}$, $B\in\Re^{N_t\times n_u}$, 
$G_d\in\Re^{N_t\times n_d}$, $E\in\Re^{N_s\times n_u}$ and 
$E_d\in\Re^{N_s\times n_d}$}.

The constraints on $x_k$ and $u_k$ can be concisely written as
\begin{subequations}
\begin{align}
x_{\min} \leq & x_k \leq x_{\max},\label{eq:state-constraints}\\
u_{\min} \leq & u_k \leq u_{\max},\label{eq:input-constraints}
\end{align}
\end{subequations}
{ with $x_{\min},\, x_{\max}\in\Re^{N_t}$ and $u_{\min},\, u_{\max}\in\Re^{n_u}$
and $\leq$ is meant in the element-wise sense.}
These constraints encompass~\eqref{eq:level-constraints} 
and~\eqref{eq:flow-constraints}.

\subsection{Demands and Electricity Prices}\label{sec:demand-electricity}
Water demand has been the main source of uncertainty for the operation
of drinking water networks and a lot of attention has been paid on
the development of models for its prediction.
Prediction methodologies span from simple linear models~\citep{Pedregal20071710}
to neural networks~\citep{Ghiassi2008,Romano2014265} and support vector
machines~\citep{PeaGuzmn2016,MsiNelMar07},
nonlinear multiple linear regression~\citep{Yasar2012},
Holt-Winters-type models~\citep{Wang+2016},
as well as more complex neuro-fuzzy models~\citep{PapPocLas16,Lertpalangsunti199921}.
Increased predictive ability can be obtained using exogenous information such as
weather forecasts~\citep{BakVanDui+14c} and calendar data~\citep{SamGroSop+14}.

In the context of a deregulated wholesale energy market, prices on the day-ahead
market are volatile and are often decided on the basis of an auction (bid-based market)
among energy companies instead of bilateral agreements with an energy provider.
In such cases, energy prices may change on a daily or hourly basis~\citep{PugPatBerBem13}.
It is then necessary to be able to predict the day-ahead evolution of the 
prices using past data; several time series analysis methodologies have 
been developed for that purpose --- see~\citep{Weron20141030} and 
references therein.

In our approach the prediction procedure is decoupled from the control
{
which allows the use of any forecasting methodology without 
having to modify the controller parameters or implementation.}
An independent \textit{forecaster} provides estimates of future demands and
electricity prices along with an estimation of their uncertainty which
is discussed in the next section. At time $k$, the predicted demands
for the future time $k+j$ are estimated by a model which computes
$\hat{d}_{k+j {}\mid{} k}$. Likewise, we denote the predicted electricity
prices by $\hat{\alpha}_{k+j {}\mid{} k}$. Let $d_{k+j}$ and $\alpha_{k+j}$
denote the actual, but unknown at time $k$, values of the water demands
and prices. Then
\begin{align}
 \label{eq:prediction-error}
 \begin{bmatrix}
  d_{k+j}\\
  \alpha_{k+j}
 \end{bmatrix}
 =
 \begin{bmatrix}
  \hat{d}_{k+j {}\mid{} k}\\
  \hat{\alpha}_{k+j {}\mid k}
 \end{bmatrix}
 +
 \epsilon_{j{}\mid{} k},
\end{align}
where $\epsilon_{j{}\mid{} k}$ is a random variable which corresponds to the $j$-step-ahead
prediction error at time $k$.
At time $k$, a forecaster provides finite-horizon estimates of the upcoming
water demands and electricity prices
\begin{align}
 \bm{\hat{d}}_k = (
  \hat{d}_{k+1 {}\mid{} k},
  \hat{\alpha}_{k+1 {}\mid k},
  \ldots,
  \hat{d}_{k+H_p {}\mid{} k},
  \hat{\alpha}_{k+H_p {}\mid k}),
\end{align}
where $H_p$ is a \textit{prediction horizon}.
This information will then be provided to the controller as we shall
discuss in Section~\ref{sec:smpc}.

\subsection{Uncertainty}\label{sec:uncertainty}
% Motivation: why do we use the scenario tree approach?
It is common in stochastic control-oriented modeling to assume that the
errors $\epsilon_{j{}\mid{}k}$ are independently distributed~\citep{WanOcaPui16,GroVelOca+16}.
This assumption however neglects the covariance across the times stages ---
indeed, if at the future time $j=1$ the model has a large prediction error
we would rather expect that the prediction error at time $j=2$ is likely to
be large too. This motivates the use of scenario trees: discrete representations
of the random processes $(\epsilon_{j{}\mid{}k})_j$ which capture such multistage
covariances~\citep{ShaDen09}.

To date, stochastic modeling for drinking water networks in presence of price
uncertainty has received little attention --- to the best of the authors'
knowledge,~\citep{Eck2014} is the only relevant reference --- and the scenario tree
approach has not been used previously.

\begin{figure}[!t]
 \centering
 \includegraphics[width=0.65\columnwidth]{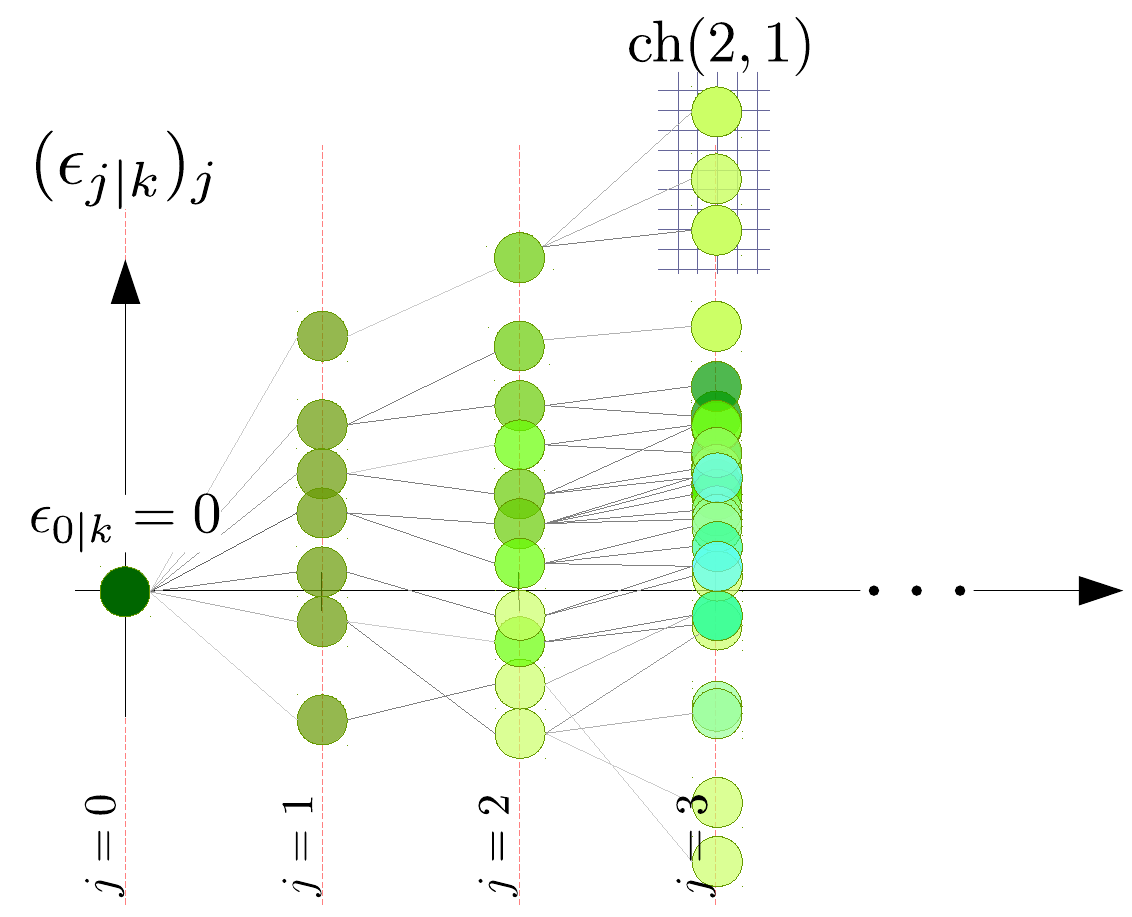}
 \caption{Scenario tree of the random process $(\epsilon_{j{}\mid{}k})_j$.
          At time $j=0$, the prediction error is $\epsilon_{0\mid k}=0$; this defines
          the \textit{root node} of the tree. {
          The children nodes $\child(2,1)$ of
          node $(2,1)$ are highlighted in the figure.}}
 \label{fig:tree}
\end{figure}

%What is a scenario tree?
A scenario tree is a structure such as the one illustrated in Fig.~\ref{fig:tree}.
The scenario tree is organised into time instants $j=0,\ldots,H_p$ called \textit{stages} and a number
of \textit{nodes} at each stage denoted by $\epsilon_{j{}\mid{}k}^i$ --- these are treated as the
possible values of $\epsilon_{j{\mid}k}$. At time $j=0$, the prediction error is always equal to
$0$ assuming that we observe the current water demands and electricity prices. The corresponding
unique node is called the \textit{root} node of the tree. The nodes of the tree at the last
stage are known as \textit{leaf nodes}. The number of nodes at stage $j$ is denoted by $\mu_j$.
All nodes are identified by a pair $\alpha=(j,i)$, where $j=0,\ldots, H_p$ is the stage index and
$i$ is the index of the node in that stage.
Each non-leaf node defines a (nonempty) set of \textit{children}
$\child(j,i)$ which are those nodes in next $j+1$ which are linked to $(j,i)$.
Conversely, each not except for the root node defines a unique \textit{ancestor}
denoted by $\anc(j,i)$. The probability of visiting a node $(j,i)$ starting from
the root node and following the tree structure is denoted by $p_j^i$.

Note that the joint demand-price modeling of the uncertainty allows us to
cast possible demand-price correlations as in cases of uncertain volume-based pricing.

Using the scenario tree structure, equation~\eqref{eq:prediction-error} yields
\begin{align}
\begin{bmatrix}
  d_{k+j}^i\\
  \alpha_{k+j}^i
 \end{bmatrix}
 =
 \begin{bmatrix}
  \hat{d}_{k+j {}\mid{} k}\\
  \hat{\alpha}_{k+j {}\mid k}
 \end{bmatrix}
 +
 \epsilon_{j{}\mid{} k}^i,
\end{align}
{ where $j=0,\ldots, H_p$ and $i=1,\ldots,\mu_j$. 
  Here we see that the tree structure of $\epsilon_{j\mid k}^i$
  induces a corresponding tree structure upon the water demands
  and electricity prices, namely $d_{k+j}^i$ and $\alpha_{k+j}^i$.
  These are the contingent future water demand and electricity price
  values associated with the prediction error $\epsilon_{j\mid k}^i$.
}

Similarly, equation~\eqref{eq:input-disturbance-coupling} gives
\begin{align}
 Eu_{k+j{}\mid{}k}^i + E_d d_{k+j{}\mid{}k}^i = 0,
\end{align}
where $(u_{k+j{}\mid{}k}^i)_{j,i}$ become the decision variables of the
stochastic optimal control problem we shall present in the following section.
The discrete-time system dynamics~\eqref{eq:dynamics} becomes
\begin{align}
 \label{eq:dynamics:tree}
 x_{k+j+1{}\mid{}k}^{l} = A x_{k+j{}\mid{}k}^{i} + Bu_{k+j{}\mid{}k}^{l} + G_d d_{k+j{}\mid{}k}^{l},
\end{align}
where $(j+1,l)\in\child(j,i)$.

Scenario trees can be constructed from observed sequences of prediction errors
which can be easily obtained in practice using methodologies such as~\citep{Pflug2015scen}
or the popular \textit{scenario reduction} method~\citep{HeiRom09}.
There exist several other scenario generation algorithms such as
clustering-based algorithms~\citep{Latorre20071339,Chen:2014:KST:2787349.2787351}
and simulation and optimization-based approaches~\citep{Beraldi20102322,Gulpınar20041291}.
It is not necessary, however, to update the scenario tree for $\epsilon_{j{}\mid{}k}$
at every time instant $k$ --- it should be updated occasionally to
detect changes in the predictive ability of the forecaster or
whenever the predictive model is updated.

\section{Scenario-based stochastic optimal control}\label{sec:smpc}
\subsection{Control Objectives}
The cost for the operation of the water network is quantified in
terms of three individual costs which have been proposed in the
literature~\citep{SamGroSop+14,ConPuiCem14,OcPuCe+09,martinez2015transport}:
the economic cost which is related
to the treatment cost and electricity required for pumping, the
smooth operating cost which penalizes the abrupt operation of
pumps and valves and the safety storage cost which penalizes the
use of water from the reserves (i.e., allowing the level in
the tanks to drop below the safety level).

The \textit{economic cost} quantifies the \textit{production} and
 \textit{distribution} cost and it is computed by
 \begin{align}\label{cost:linear}
  \ell^{w}(u_k,k)=W_{\alpha}(\alpha_0+\alpha_{k})'u_k,
 \end{align}
where $\alpha_0' u_k$ is the water production cost (treatment and acquisition fees),
$\alpha_{k}' u_k$ is the uncertain pumping cost and
$W_{\alpha}$ is a positive scaling factor.

The \textit{smooth operation cost} is defined as
\begin{align}\label{cost:smooth}
 \ell^{\Delta}(\Delta u_k)=\Delta u_k' W_u \Delta u_k,
\end{align}
where $\Delta u_k=u_k-u_{k-1}$ and $W_u\in\Re^{n_u\times n_u}$ is a symmetric
positive definite weight matrix.

The total \textit{stage cost} at a time instant $k$ is the summation of the above
costs and is given by
\begin{align*}
 \ell(u_k, u_{k-1}, k)=\ell^{w}(u_k,k)+\ell^{\Delta}(\Delta u_k).
\end{align*}

The \textit{safety storage cost} penalizes the drop of water level
in the tanks below a given \textit{safety level} $x_s$. An elevation
above this safety level ensures that there will be enough water
in unforeseen cases of unexpectedly high demand
and also maintains a minimum pressure for the flow of water
in the network. This is given by
\begin{align}\label{cost:non-smooth}
 \ell^{S}(x_k)= W_s \|\max\{0, x_s-x_k\}\|,
\end{align}
where $W_s$ is a positive scaling factor.

The state constraints~\eqref{eq:state-constraints} should be satisfied at
all times without however jeopardizing the feasibility of the optimal
control problem we have to solve at every time instant. For that,
we introduce an additional cost which penalizes the violation of
the state constraints as follows
\begin{align}
 \ell^x(x_k) &= W_x\big(\|\max\{0,x_{\min}-x_k\}\|\notag\\
 &\qquad\ + \|\max\{0, x_k-x_{\max}\}\|\big),
\end{align}
where $W_x$ is a positive weight factor.

The scaling factors $W_{\alpha}$, $W_u$, $W_s$ and $W_x$ are the
tuning knobs of stochastic MPC as we shall discuss in the following
section.

\subsection{Stochastic optimal control: problem formulation}
\begin{figure}[t]
 \centering
 \includegraphics[width=0.6\columnwidth]{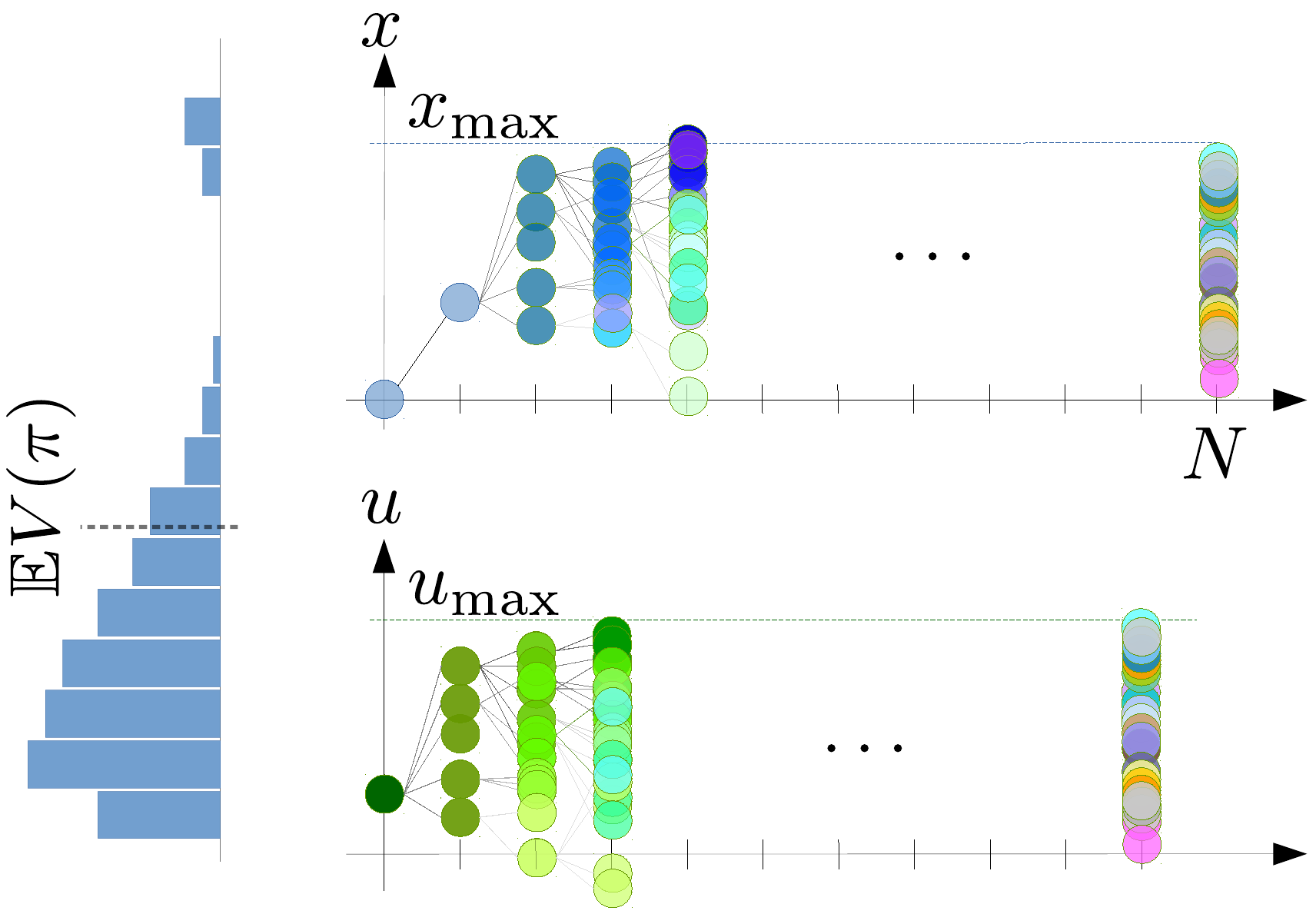}
 \caption{The concept of scenario-based stochastic optimal control:
          at every time instant $k$, we make an optimal contingency
          plan by minimizing the \textit{expectation} of a
          cost function $V$ which encodes the operation cost along a
          finite prediction horizon.}
 \label{fig:smpc}
\end{figure}

In scenario-based stochastic MPC, at every time instant $k$ we solve a
stochastic optimal control problem which consists in determining an optimal
\textit{contingency plan} for the future course of actions in a \textit{causal
fashion}, that is, our future decisions $u_{k+j{}\mid{}k}$ are only allowed
to depend on information that will be available to the controller at time
$k+j$~\citep{BerShr96}. This was tacitly stated in equation~\eqref{eq:dynamics:tree}.
This concept is illustrated in Fig.~\ref{fig:smpc}.

We formulate the following scenario-based stochastic MPC problem
with a prediction horizon $H_p$ and decision variables
$\bm{x}=\{u_{k+j \mid k},x_{k+j+1 \mid k}\}_{j=0,\ldots,H_p-1}$
where we minimize the expected total cost along the prediction horizon
\begin{subequations}\label{eq:SMPC}
\begin{align}
 &\minimise_{\bm{x}}\ \E \left[ V(\bm{x},p,q,k) \right] + V_s(\bm{x}),\label{eq:smpc-cost}
\end{align}
subject to the following constraints
\begin{align}
& x_{k\mid k}=p,\ u_{k-1\mid k}=q,\label{eq:inital_SP}\\
& x_{k+j+1{}\mid{}k}^{l} = A x_{k+j{}\mid{}k}^{i}
    + Bu_{k+j{}\mid{}k}^{l} + G_d d_{k+j{}\mid{}k}^{l},\label{eq:SMPC:dynamics}\\
& Eu_{k+j \mid k}^{i}+E_dd_{k+j\mid k}^{i}(\epsilon_j)=0,\label{eq:SMPC:mass-pres}\\
& u_{\min}\leq u_{k+j\mid k} \leq u_{\max}^{i}.
\end{align}
\end{subequations}
In~\eqref{eq:smpc-cost}, $\E$ is the expectation operator, $V$
is the cost function given by
\begin{align}
V(\bm{x},p,q,k)=\sum_{j=0}^{H_p-1}\ell(u_{k+j \mid k},u_{k+j-1 \mid k}, k{+}j),
\end{align}
and $V_s$ is the total state constraint violation penalty
defined as
\begin{align}
 V_{s}(\bm{x}) = \sum_{j=0}^{H_p}\sum_{i=1}^{\mu_{j-1}}
          \ell^x(x_{k+j\mid k}^{i}) + \ell^S(x_{k+j\mid k}^{i}).
\end{align}
Note that we have a time-varying cost as the electricity prices change
with time.

% \begin{enumerate}
%  \item Problem statement: emphasise that this an entirely data-driven approach
% \end{enumerate}

\section{Numerical Algorithm}
\subsection{Problem reformulation}
Most modern numerical optimization algorithms such as the (accelerated)
proximal gradient algorithm, the alternating directions method of multipliers
(ADMM)~\citep{ParBoy13}, the Pock-Chambolle method~\citep{chambollepock11},
Tseng's forward-backward-forward algorithm~\citep{Tseng2000} and many another
require that the optimization problem be first written in a form
\begin{align}
 \label{eq:generic-optimization-problem}
 \mathcal{P}:\ \minimise_{\bm{x}\in\Re^n}\, f(\bm{x}) + g(H\bm{x}),
\end{align}
where $f:\Re^n\to\barre \dfn \Re\cup\{+\infty\}$ and $g:\Re^m\to\barre$ are convex,
lower semi-continuous extended-real-valued functions and $H:\Re^n\to\Re^m$ is a linear
operator. Functions $f$ and $g$ are allowed to return the value $+\infty$ to encode
constraints; for example, the constraint $x\in C$ is encoded by the 
\textit{indicator function} of the set $C$ which is
\begin{align}
 \label{eq:indicator-function}
 \delta(x{}\mid{}C) = \begin{cases}
                       0,&\text{ if } x\in C\\
                       +\infty,&\text{ otherwise}
                      \end{cases}
\end{align}

The key question is how to split the optimization problem in~\eqref{eq:SMPC}
so that the resulting formulation is amenable to a fast numerical solution
with massive parallelization. For reasons that will be elucidated in
Section~\ref{sec:apg-algorithm} we choose $f$ to be the \textit{smooth}
part of the cost (which corresponds to the linear function $\ell^w$ and
the quadratic function $\ell^\Delta$) plus the indicator of the input-disturbance
coupling given in~\eqref{eq:input-disturbance-coupling} plus the indicator
of the system dynamics in~\eqref{eq:dynamics}, that is $f:\Re^n\to\barre$ is
defined as
\begin{align}
 f(\bm{x}) &= \delta(u_j^i|\Phi_1(d_j^i))
	    + \delta(x_{j+1}^{i}, u_j^i, x_j^{\anc(j+1,i)}|\Phi_2(d_j^i))\notag\\
	   &+\sum_{j=0}^{H_p{-}1}\sum_{i=1}^{\mu_{j}}
	      p_j^i(\ell^w(u_j^i)+\ell^{\Delta}(\Delta u_j^i)),\label{eq:f}
\end{align}
where $\Delta u_j^i= u_j^i- u_{j-1}^{\anc(j,i)}$ and
$\Phi_1(d)$ is the affine subspace of $\Re^{n_u}$
\begin{equation}\label{eq:Phi_1}
\Phi_1(d) = \{u:Eu + E_d d = 0\},
\end{equation}
and $\Phi_2(d)$ is the affine subspace of $\Re^{2n_x+n_u}$
defined by the system dynamics
\begin{equation}\label{eq:Phi_2}
\Phi_2(d)=\{(z, x, u): z=Ax + Bu + G_d d\}.
\end{equation}

Function $g$ is naturally chosen to be the indicator of the set of
input constraints plus the total constraint violation penalty function $V_s$.
Note, however, that the same variable $x_{k+j{}\mid{}k}^{i}$ participates
in both functions $\ell^x$ and $\ell^S$.
As we shall explain in Section~\ref{sec:apg-algorithm}, this complicates any
computations thereon. For that reason we introduce a linear operator $H:\bm{x}\mapsto \bm{y} \dfn H(\bm{x})$
which maps $x_{k+j{}\mid{}k}^i$ to $(x_{k+j{}\mid{}k}^i, x_{k+j{}\mid{}k}^i)$
and $u_{k+j{}\mid{}k}^i$ to itself, that is
\begin{subequations}
\label{eq:H}
\begin{align}
 y_{k+j{}\mid{}k}^i = (x_{k+j{}\mid{}k}^i, x_{k+j{}\mid{}k}^i, u_{k+j{}\mid{}k}^i),
\end{align}
for $j=0,\ldots,H_p-1$ and
\begin{align}
 y_{k+H_p{}\mid{}k}^i = (x_{k+H_p{}\mid{}k}^i, x_{k+H_p{}\mid{}k}^i).
\end{align}
\end{subequations}
Then, we define the function $g:\Re^m\to\barre$
\begin{align}
 \label{eq:g}
 g(\bm{y}) &= \sum_{
 \substack{
 j=0,\ldots,H_p\\
 i=1,\ldots,\mu_{j-1}}}
          \ell^x(y_{k+j\mid k}^{i,(1)}) + \ell^S(y_{k+j\mid k}^{i,(2)})\notag\\
          &+ \sum_{
 \substack{
 j=0,\ldots,H_p-1\\
 i=1,\ldots,\mu_{j-1}}}\delta(y_{k+j\mid k}^{i,(3)} {}\mid{} \mathcal{U}),
\end{align}
where $\mathcal{U} \dfn \{u\in\Re^{n_u} {}\mid{} u_{\min} \leq u \leq u_{\max}\}$.

The scenario-based optimization problem~\eqref{eq:SMPC} is now
in the form~\eqref{eq:generic-optimization-problem} with $f$ and $g$
given by~\eqref{eq:f} and~\eqref{eq:g} respectively and $H$ given
by~\eqref{eq:H}.

\subsection{Convex conjugates and proximal operators}
Before we can proceed with the statement of the numerical algorithm for the
solution of problem~\eqref{eq:generic-optimization-problem} we need to
introduce a few mathematical notions. A function $f:\Re^n\to\barre$
is called \textit{proper} if it is not everywhere equal to $+\infty$. It is called
\textit{lower semi-continuous} if for every $x\in\Re^n$, $\liminf_{z\to x} f(z) = f(x)$.
The \textit{domain} of $f$ is the set $\dom f \dfn \{x\in\Re^n {}\mid{} f(x) < \infty\}$.
It is called $\kappa$-strongly convex, for some $\kappa>0$, if the function
$f(x)-\nicefrac{\kappa}{2}\|x\|^2$ is convex.

For a proper, convex, lower semi-continuous function $f$, we define its
\textit{convex conjugate} to be the convex function $f^{*}:\Re^n\to\barre$
defined as~\citep[Ch.~11]{RocWets09}
\begin{align}
 f^*(y) = \sup_{x\in\Re^n}\  y'x - f(x).
\end{align}

An important property is that if $f$ is $\kappa$-strongly convex, then
$f^{*}$ is differentiable with $\nicefrac{1}{\kappa}$-Lipschitz gradient%
~\citep[Prop.~12.60]{RocWets09} and the gradient of $f^{*}$ is
\begin{align}
 \label{eq:dual-gradient}
 \nabla f^{*}(y) = \argmax_{x\in\Re^n}\  y'x - f(x) .
\end{align}

By means of the convex conjugate we may derive the \textit{Fenchel dual}
optimization problem of~\eqref{eq:generic-optimization-problem} which
is~\citep[Sec.~15.3]{BauCom2011}
\begin{align}
 \label{eq:fenchel-dual}
 \mathcal{D}:\ \minimise_{y\in\Re^m} f^*(-H'y) + g^*(y).
\end{align}

Fenchel duality offers a more powerful framework as compared to the
classical Lagrangian duality approach as it allows us to dualise
\textit{functions} (by means of their convex conjugates) rather than
merely constraints.

Under certain conditions on $f$ and $g$ (see Section~\ref{sec:apg-algorithm}), the two problems are equivalent:
their optimal values are equal and given a dual-optimal point $y^\star$
which is a minimizer of~\eqref{eq:fenchel-dual}, the optimal solution of $\mathcal{P}$
is~$x^\star = \nabla f^{*}(-H'y^\star)$.
This is referred to as \textit{strong duality}.
The reason why we formulate the dual optimization problem $\mathcal{D}$
is because it possesses a favorable structure which can be exploited
in the development of fast and parallelizable numerical algorithms
(See Section~\ref{sec:apg-algorithm}).

Lastly, for a proper, lower semi-continuous, extended-real valued function $g:\Re^m\to\barre$
we define its \textit{proximal operator} with parameter $\gamma>0$ to be a
function $\prox_{\gamma g}:\Re^m\to\Re^m$ defined as
\begin{align}
 \prox_{\gamma g}(v) \dfn \argmin_{x\in \Re^m} g(x) + \nicefrac{1}{2\gamma}\|x-v\|^2.
\end{align}
Proximal mappings act as generalized projections. For example, the proximal mapping
of the indicator function $\delta({\,}{\cdot}{\,}{\mid}{\,}C)$ --- cf.~\eqref{eq:indicator-function} ---
of a nonempty, closed, convex set is the projection onto that set, i.e.,
$\prox_{\gamma \delta({}{\cdot}{}\mid{} D)}(v) = \proj(v{}{\mid}{}C)$.

The proximal operators of many convex functions (such as Euclidean norm, norm-1, quadratic
and linear functions, distance-to-set functions) are easy to compute and typically
consist in element-wise operations which can be fully parallelized on a GPU.
We shall refer to such functions as \textit{prox-friendly}~\cite{Wytock2016}.

So long as $\prox_{\gamma g}$ is easy to compute, so is $\prox_{\gamma g^*}$
and it can be obtained from the \textit{Moreau decomposition formula} which is
\begin{align}
 \label{eq:moreau-decomposition}
 \prox_{\gamma g}(v) + \gamma \prox_{\gamma^{-1} g^*}(\nicefrac{v}{\gamma}) = v.
\end{align}

\subsection{Accelerated proximal algorithm on the dual optimization problem}
\label{sec:apg-algorithm}
The function $f(\bm{x}) + g(H\bm{x})$, with $f$ and $g$ defined by~\eqref{eq:f} and~\eqref{eq:g}
respectively, is proper, convex and piecewise linear-quadratic on its domain, so,
following~\cite[Thm.~11.42]{RocWets09} there is strong duality.
As discussed above, since $f$ is strongly convex, $f^*$ is differentiable with Lipschitz gradient.
Function $g$ is written in the form of a \textit{separable sum} ---  a sum
of prox-friendly functions of different arguments~\cite{ParBoy13}.
Function $g$ is indeed prox-friendly. Let $(\prox_{\gamma g}(v))_{j,i,(s)}$ denote
the part of the vector $\prox_{\gamma g}(v)\in\Re^m$, indexed by $j$, $i$ and $s$,
which corresponds to $y_{k+j{}\mid{}k}^{i,(s)}$, for $s=1,2,3$. Then, $(\prox_{\gamma g}(v))_{j,i,(1)}$
and $(\prox_{\gamma g}(v))_{j,i,(2)}$ are computed by virtue of the formula
\begin{align}
 &\prox_{\gamma\dist(\cdot\mid C)}(v)  \notag\\
 &=
    \begin{cases}
      v
	  + \frac{\proj_{C}(v)-v}
		{\dist(v {}\mid{} C)},&\text{if } \dist(v{}{\mid}{}C) > \gamma\\
    \proj_{C}(v),&\text{otherwise}
      \end{cases}
\end{align}
where $\dist(\cdot{}\mid{}C)$ is the distance-to-set function and
the fact that $\ell^x(x) = W_x\dist(x \mid [x_{\min},x_{\max}])$ and
$\ell^S(x) = W_s \dist(x\mid [x_s, +\infty))$.
The proximal operator $(\prox_{\gamma g}(v))_{j,i,(3)}$ is simply the
projection on $\mathcal{U}$.

Note that although $g$ is prox-friendly, $g$ composed with the linear operator $H$
--- as it is in~\eqref{eq:generic-optimization-problem} --- is not.
This is the main reason why we resort to the dual problem~\eqref{eq:fenchel-dual}.

Given the properties of functions $f^*$, being differentiable with Lipschitz
gradient, and $g^*$, being prox-friendly, we may use Nesterov's accelerated
proximal gradient method on the dual problem which produces the sequence
\begin{subequations}
\begin{align}
 w^{\nu} &= y^{\nu} + \beta_{\nu}(y^\nu - y^{\nu-1}),\label{eq:apg:1}\\
 y^{\nu} &= \prox_{\gamma g^*}(w^\nu + \gamma H\nabla f^*(-H'w^\nu)),\label{eq:apg:2}
\end{align}
\end{subequations}
with $y^0=0$, $y^{-1}=0$, $\beta^0=0$.
In~\eqref{eq:apg:1} we perform an \textit{extrapolation step} for some $\beta_\nu>0$
and we perform a dual gradient projection update on the extrapolated vector $w^\nu$.
The extrapolation parameters $\beta-\nu$ are parametrized as $\beta^\nu=\theta_\nu(\theta_{\nu-1}^{-1}-1)$
with $\theta_0=\theta_{-1}=1$. Any choice of $\theta_\nu$ so that
$1-\theta_{\nu+1} \leq \nicefrac{\theta_{\nu+1}^2}{\theta_{\nu}^2}$.
A simple choice is
\begin{align}
 \theta_{\nu} = \frac{2}{\nu+2},
\end{align}
for $\nu\geq 1$. Here we choose $\theta_{\nu+1} = \nicefrac{1}{2}(\sqrt{\theta_\nu^4 + 4 \theta_\nu^2}-\theta_\nu^2)$
which satisfies the above requirement with equality.

This algorithm fits into Tseng's Alternating Minimization framework~\citep{Tseng1991,Tseng1990}.
It has been found to be suitable for embedded applications as it is relatively
simple to implement and it has good convergence properties (the dual variable $y^{\nu}$ converges
with rate $O(\nicefrac{1}{\nu^2})$ and an averaged primal iterate converges at $O(\nicefrac{1}{\nu^2})$
as well)~\citep{PatBem2014}.
It involves only matrix-vector operations (additions and multiplications)
and it is numerically stable. In the following section we discuss how
the involved operations can be massively parallelized in a lock-step fashion
(performing the exact same operation on different memory positions)
and how the algorithm can be implemented on a GPU.

\subsection{Implementation}
Because of the definition of $f$ in~\eqref{eq:f}, the computation of the gradient of
$f^*$ as defined in problem~\eqref{eq:dual-gradient} boils down to the solution of
a scenario-based \textit{optimal control} problem where the only constraints are
the ones defined by the system dynamics. This problem can be solved by dynamic programming
leading to a Riccati-type recursion from stage $k=H_p$ to stage $k=0$.
At each stage operations across all nodes can be fully parallelized.
In particular, such a parallelization --- assuming that full parallelization is
supported by the hardware --- equalizes the complexity of the scenario-based
Riccati recursion to that of a deterministic one.
A detailed exposition of the details of this procedure
is available in~\cite{SamSopBemPat15,SamSopBemPat17b}.

All operations involved in the computation of $\prox_{\gamma g}$ are element-wise
operations and can be fully parallelized on a GPU; therefore, the
computational cost for applying $\prox_{\gamma g}$ --- or, what is the same ---
$\prox_{\gamma^{-1} g^*}$ via~\eqref{eq:moreau-decomposition} is
negligible.

\section{Functionality}
In this section we present \texttt{\name}, a \texttt{CUDA-C++} implementation of
the accelerated proximal gradient method for the solution of scenario-based
stochastic optimal control problems, particularly tailored to the needs of a
drinking water network.

\subsection{Software structure}
\begin{figure}
 \centering
 \includegraphics[width=0.7\columnwidth]{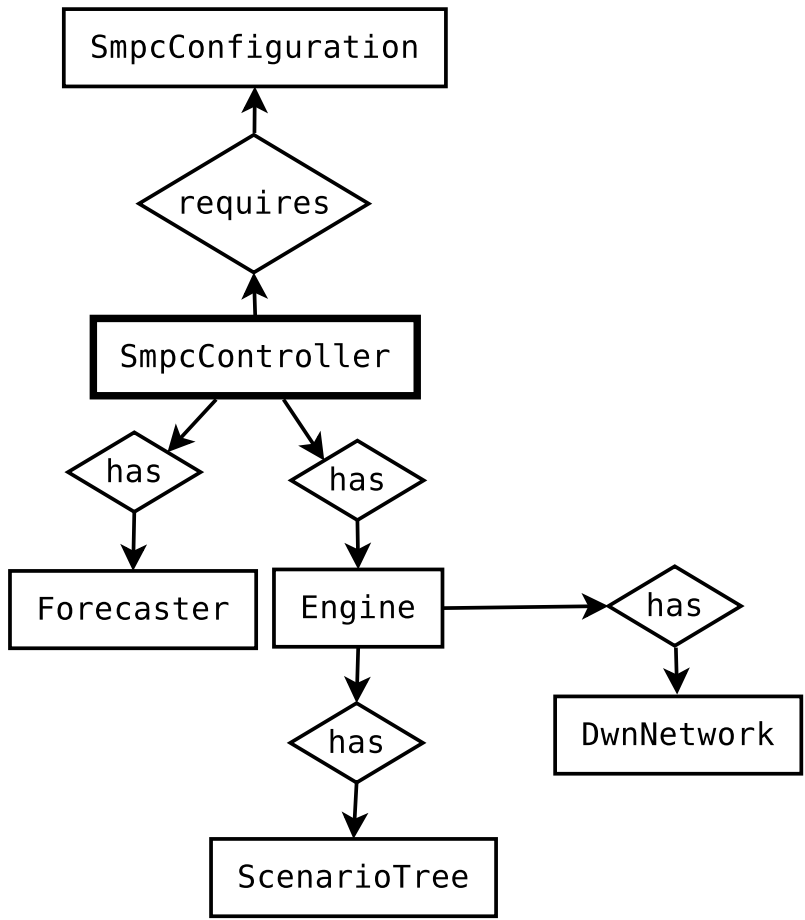}
 \caption{{Entity-Relationship (ER) diagram of the} entities in \texttt{\name} which reflects the underlying class structure.}
 \label{fig:er}
\end{figure}
The entities involved in \texttt{\name} and the relationships
among each other are illustrated in Fig.~\ref{fig:er} which
correspond to classes in the \texttt{CUDA-C++} implementation of
\texttt{\name} (see also Table~\ref{table:Class}).

\begin{table*}[!t]
\centering
\footnotesize
 \caption{Description of the major class and the associated JSON files
          in \texttt{\name}.}\label{table:Class}
 \begin{tabular}{p{2.5cm}|p{5cm}|p{3.3cm}}
 \textbf{Class name} & \textbf{Description} & \textbf{JSON file} \\
  \hline
 \texttt{SmpcController} & 
		Runs the accelerated dual proximal gradient algorithm can computes
                a control action to be applied to the water network. 
                Computations are carried out on GPU and the output (flow 
                set-points) can be stored in a JSON file.
                & \texttt{controlOutput.json}\\
 \hline
 \texttt{DwnNetwork}  & Encapsulates all information related to the topology, 
			dynamics and constraints of the water network.
                      & \texttt{network.json}\\
 \hline
 \texttt{ScenarioTree}  & Scenario-tree representation of the uncertainty
                          in electricity prices and water demands.
                        & \texttt{scenarioTree.json} \\
 \hline
 \texttt{Forecaster} & An abstract forecaster which predicts the upcoming electricity
  		       prices and water demands using some predictive model (implemented by
 		       subclassing \texttt{Forecaster}) or reads the forecasts from a
 		       JSON file (so that the user can use forecasts from 
 		       third-party software).
		     & \texttt{forecaster.json}\\
 \hline
 \texttt{Engine}  &  Provides essential functionality to \texttt{SmpcController} 
                     and manages the GPU-side memory. &\\
 \hline
 \texttt{SmpcConfiguration} & Contains configuration parameters that are relevant for
                     \texttt{SmpcController} (tuning parameters, solver tolerance, maximum number of iterations). & \texttt{controllerconfig.json}\\
 \end{tabular}\\[1.5em]
 \hrule
\end{table*}

The \texttt{DwnNetwork} class stores data related to the
network topology, physical constraints and network dynamics matrices
in equations~\eqref{eq:dynamics} and~\eqref{eq:input-disturbance-coupling}.
The user can create instances of \texttt{DwnNetwork} simply by passing a JSON file,
\texttt{network.json}, with the network information.

The \texttt{ScenarioTree} class models the scenario tree structure that
represents the uncertainty associated with the volatile energy prices and
water demands.
This class describes the structure of the scenario tree by assigning a
unique index to each node and storing the indexes of the children of 
each non-leaf node, the ancestor of all nodes except for the root 
node and the values $d_{k+j\mid k}^i$ and $\alpha_{k+j\mid k}^i$
at each node $(j,i)$.
An instance of \texttt{ScenarioTree} can be generated using a JSON
file, \texttt{scenarioTree.json}.

An \texttt{SmpcController} executes the accelerated proximal
gradient algorithm to solve a stochastic optimal control problem
at every time instant and compute the control actions to be applied
to the water network. 
{ Certain functionality is delegated to
\texttt{Engine} --- a collection of utility methods --- 
which precomputes certain quantities which
are associated with the Riccati-type recursion (details can be 
found in~\cite{SamSopBemPat17b}) for the computation of
the dual gradient and manages the associated memory on GPU.
}

The Engine is, in turn, linked with to a \texttt{DwnNetwork} entity which
provides all necessary technical specifications for the network
(topology, dynamics, constraints)
and a \texttt{ScenarioTree}
entity which encodes the probabilistic information associated
with the prediction errors.
{ Note that the end-user does not have to create instances
of \texttt{Engine} or directly interact with it.}

A \texttt{Forecaster} provides to the controller estimates of
the upcoming water demands and electricity prices.
This is an abstract class which can be subclassed with particular
model implementations (e.g., ARIMA, SVM, or any other), or
the user can provide custom forecasts using any third-party software
which exports its forecasts in a JSON file.

\texttt{SmpcController} requires certain configuration
parameters which is provided by the entity \texttt{SmpcConfiguration}.
There, the user specifies the desired tolerance, maximum number of
iterations and can override other solver-specific properties.

\begin{figure}
 \centering
 \includegraphics[width=0.7\linewidth]{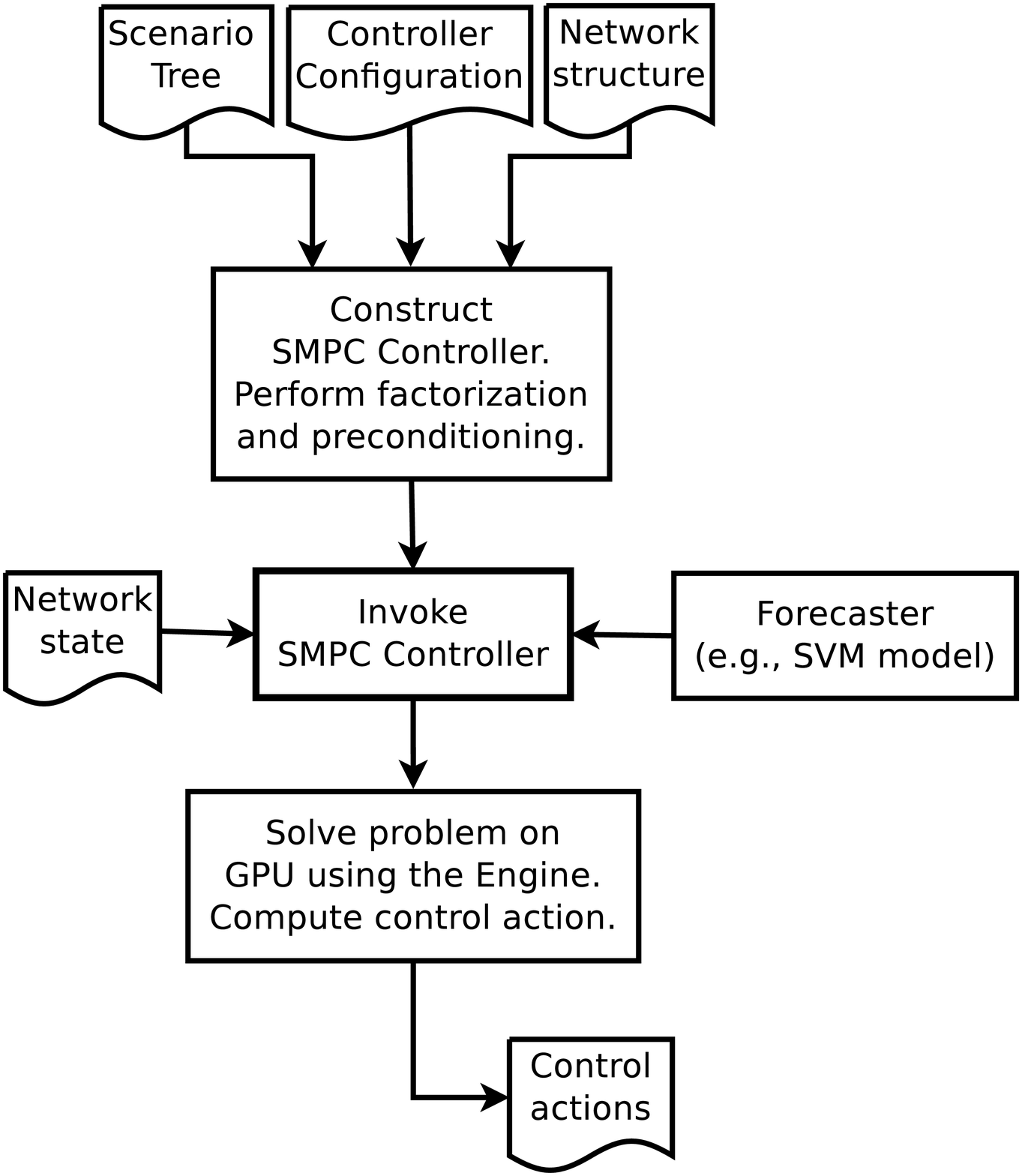}
 \caption{Flow of information in \texttt{\name}.}
 \label{fig:flow-of-info}
\end{figure}

Overall, the flow of information in \texttt{\name} is shown in
Fig.~\ref{fig:flow-of-info}. The end-user initializes an \texttt{SmpcController}
object by providing the network topology, the controller configuration 
and a scenario tree. During real-time operation, the controller receives the 
network state $x_k$ (which can be provided in a JSON file) and, using
a demand/price forecaster, computes a control action which is applied to the system.

In Table~\ref{table:functions} we list the main methods of \texttt{\name}.
Additional getter methods are available in each class for more advanced
use-case scenarios.

\begin{table*}
\centering
\footnotesize
  \caption{Key methods in \texttt{RapidNet}.}\label{table:functions}
  \begin{tabular}{p{2.5cm}|p{2.5cm}| p{5.5cm}}
  \textbf{Class} & \textbf{Method} &  \textbf{Description} \\
  \hline
  \texttt{SmpcController} & \texttt{controlAction}  & Computes control action using the accelerated proximal gradient algorithm to solve the scenario-based stochastic optimal control problem.\\
  \hline
  \texttt{Forecaster}     & \texttt{predictDemand}  & Returns water demand forecasts.\\
  \hline
                          & \texttt{predictPrice}   & Returns energy price forecasts.\\
  \hline
  \texttt{Engine}& \texttt{factorStep}  & Precomputes certain quantities that facilitate and accelerate the computation of the dual gradient in the algorithm.
  \end{tabular}\\[0.5em]
   \hrule
\end{table*}

\subsection{GPU Implementation}
GPUs were first developed for video applications and, due to the high demand 
in high-performance graphics, rapidly evolved to powerful hardware featuring
hundreds of computation cores. Nowadays, GPUs are used for more than video
processing and they are becoming popular for computational purposes including, 
but not limited to, environmental modeling~\cite{Xia201628,Le20151,Guidolin2016378}. 
By design they are well-suited for data-parallel lockstep applications
where the same type of operation is applied to different memory positions.
Instructions are sent to the GPU (from the CPU) in the form of \textit{compute
kernels}. GPUs offer unprecedented parallelization capabilities provided that 
the program can be parallelized in a lockstep fashion (the same operation is 
executed on different memory positions).

CUDA is a parallel computing framework and application programming interface 
for NVIDIA GPUs used for general-purpose computing. Part of the CUDA framework 
is \texttt{cuBLAS}, a parallel counterpart of the popular linear algebra library BLAS.

In \texttt{\name}, at every time instant $k$, the method \texttt{controlAction} in \texttt{SmpcController}
returns the control action that is to be applied to the water network (pump and valve
set points).
All computations involved in this method are
either summations or matrix-vector multiplications which can be parallelized across
the nodes of a stage. These multiplications are implemented using the
function \texttt{cublasSgemmBatched} of cuBLAS and vector additions are performed using the
\texttt{cublasSaxpy} of cuBLAS. Apart from the standard cuBLAS methods, we have
defined custom kernels for the summation over the set of nodes and to evaluate
projections with respect to the box constraints and the proximal operator of the
distance function from the set (that is, to compute the proximal operator of $g$
as discussed in Section~\ref{sec:apg-algorithm}.

\subsection{Software verification}
The validation of the software is done through \textit{unit testing}. A unit
represents the smallest functional part of the software and unit testing
involves verification of its functionality through predefined inputs and expected
outputs.
It assists in the debugging and maintenance of the code and facilitates
the integration of the various units reliably.
In lack of a standardized testing framework for CUDA applications, 
we developed our in-house testing framework.
Moreover, using \texttt{cudaMemCheck} we have thoroughly tested for GPU-side 
memory leaks.

\section{{\name} in action: Simulation results}
\begin{figure*}
 \centering
 \includegraphics[width=0.7\linewidth]{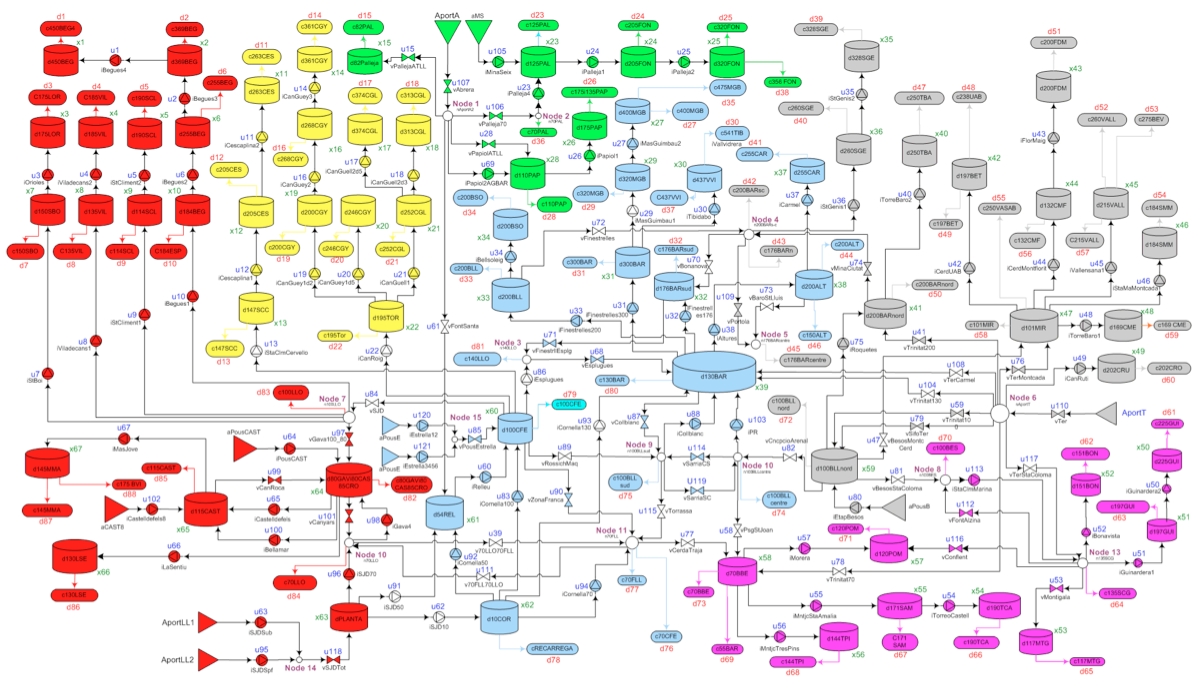}
 \caption{Schematic diagram of the drinking water network of Barcelona.
 { Overall, this distribution network involves $63$ tanks, 
 $114$ controlled flows and, in particular,
 $75$ pumps and $39$ valves, $88$ demand sectors and $17$ mixing nodes.}
          %Image source:~\citep{Rangel1}.
          }
 \label{fig:dwn}
 \hrule
\end{figure*}
In this section, we present the application of \texttt{\name} for the
management of the drinking water network of the city of Barcelona, 
whose schematic diagram is shown in Fig.~\ref{fig:dwn}, using the demand data provided
in~\cite{GroOcam+14,SamSopBemPat17b}.
The network counts a total of $63$ tanks, $114$ controlled flows 
(by means of $75$ pumps and $39$ valves), $88$ demand sectors and $17$ 
mixing nodes.

\subsection{Forecasting of water demands}
\begin{figure}[!]
 \centering
 \includegraphics[width=0.7\columnwidth]{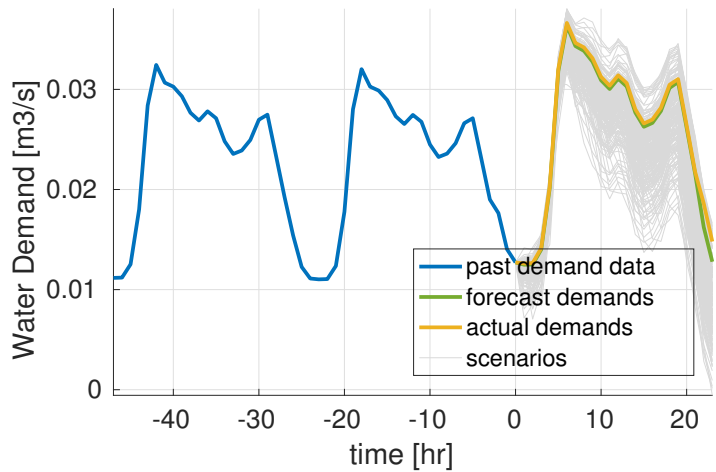}
 \caption{Water demands predicted with a radial-basis-function SVM 
 model. The current time instant is positioned at $0$.}
 \label{fig:water-demands}
\end{figure}

Upcoming water demands are predicted using a radial-basis-function 
support vector machines (SVM) model from the literature with good 
predictive ability~\cite{SamGroSop+14}.
The model predicts the water demand using past demand data together with 
calendar data (day of the week). Validation information is provided 
in~\cite{SamGroSop+14}.
In Fig.~\ref{fig:water-demands} we show an instance of a prediction using 
this SVM model.

\subsection{Forecasting of energy prices}
Among other European countries such as Denmark and Sweden, 
Austria implements a deregulated energy market which induces a volatile
stream of electricity prices. Time series of energy prices in Austria 
have become public by EXAA {
(Energy Exchange Austria; a central European energy exchange)}
and are available online at \url{http://www.exaa.at/de/marktdaten/historische-daten}
(Accessed on \today). 
{
In lack of electricity price data from Spain or relevant data from the 
water network of Barcelona, we used price data from Austria as 
an indicative dataset.
}

\begin{figure}[!]
 \centering
 \includegraphics[width=0.7\columnwidth]{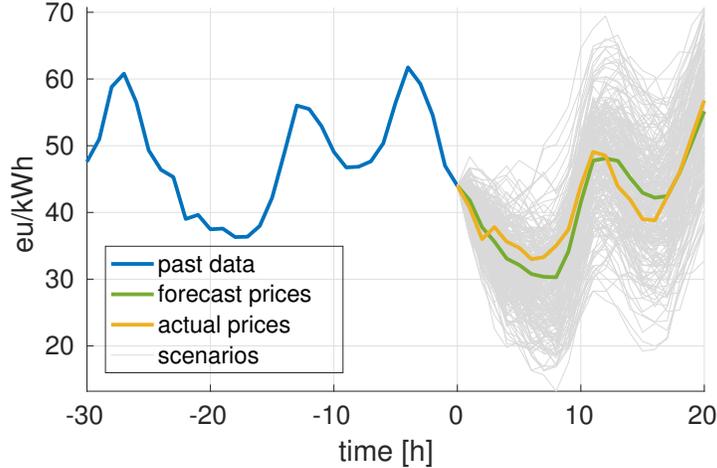}
 \caption{Past data of energy prices (blue line) together with the nominal 
          24-hour-ahead forecast produced by the ARIMA model (green line) and the actual 
          upcoming price data (orange line). A set of $250$ scenarios in 
          also shown with thin gray lines.
          The current time instant is positioned at $0$.}
 \label{fig:prices}
\end{figure}

Out of $8784$ hourly price data which are available for the year 2016, we 
excluded the last $2000$ data points to be used for testing and using the rest of 
the data we built an ARIMA$(24,1,4)$ model. 
ARIMA models have been previously used for the short-term prediction of electricity
prices in the day-ahead market~\cite{Weron20141030}.
Using a Monte-Carlo method, 
a set of $10^4$ independent scenarios were generated and, subsequently, 
these were reduced into a scenario tree using the method described in~\cite{HeiRom09}.
An instance of a prediction using the trained ARIMA model along with the associated
scenario tree is shown in Fig.~\ref{fig:prices}.

The residuals of the model where found to be uncorrelated at the  
confidence level of $99.9\%$. Indeed, the residuals pass the Ljung-Box 
Q-test of uncorrelatedness with p-value equal to $1.0000$.
The model was selected using the Akaike information criterion
(AIC) with value $1.523$.

It is natural to expect that the upcoming energy prices can be only predicted up to
moderate accuracy as they do not follow a regular pattern and are influenced
by many market-related parameters. 
Despite our limited predictive capacity, we shall show
in Section~\ref{sec:closed-loop-simulations} that by taking into account 
this volatility using stochastic
model predictive control we do mitigate the effect of the price uncertainty
leading to a more economic operation of the water network.

\subsection{Closed-loop simulations}
\label{sec:closed-loop-simulations}
In this section we present closed-loop simulation results on 
the water network where the sampling time is equal to $\unit[1]{hr}$
and the prediction horizon of the SSMPC was fixed 
to $H_p = 24$. The weight parameters used to tune the SSMPC
are $W_{\alpha} = 10^{6}$, $W_u = 1.3 \cdot 10^{7} \cdot I$, 
$W_s = 10^{5}$ and $W_x = 10^8$. Note that all units used in this 
section are SI units (flows in $\unitfrac{m^3}{s}$ and volumes in 
$\unit{m^3}$).

The system was simulated for a period of $H_s = 168$ time instants,
which corresponds to one week of operation. In order to assess the 
performance of the closed-loop operation, we use three KPIs, in particular:
(i) the \textit{economic index} which is defined as
\begin{align}
 \mathrm{KPI}_{E} = \nicefrac{1}{H_s}\sum_{k=1}^{H_s}(\alpha_0 + \alpha_k)'u_k,
\end{align}
which provides an estimation of the average hourly cost of 
operation of the water network, 
(ii) the \textit{safety index}, which is defined as
\begin{align}
 \mathrm{KPI}_{S} = \sum_{k=1}^{H_s}\| \max\{x_s - x_k, 0\} \|_1,
\end{align}
which quantifies the total weekly violation of the safety constraint (the 
``soft'' requirement that $x_k\geq x_s$), and, last, (iii) the 
\textit{complexity index} which is the worst-case time required to solve
the SSMPC optimization problem up to the specified accuracy,
which is defined as 
\begin{align}
 \mathrm{KPI}_{\tau} = \max_{k=1,\ldots, H_s}\tau_k,
\end{align}
where $\tau_k$ is the computation time required by the solver 
for solving the SSMPC problem at time instant~$k$.
The maximum runtime is of higher importance that the average
runtime in applications to verify that a decision can be 
made within the available time period.
All three indices are defined so that low values are 
preferred. 
\begin{figure}[!]
 \centering
 \includegraphics[width=0.98\columnwidth]{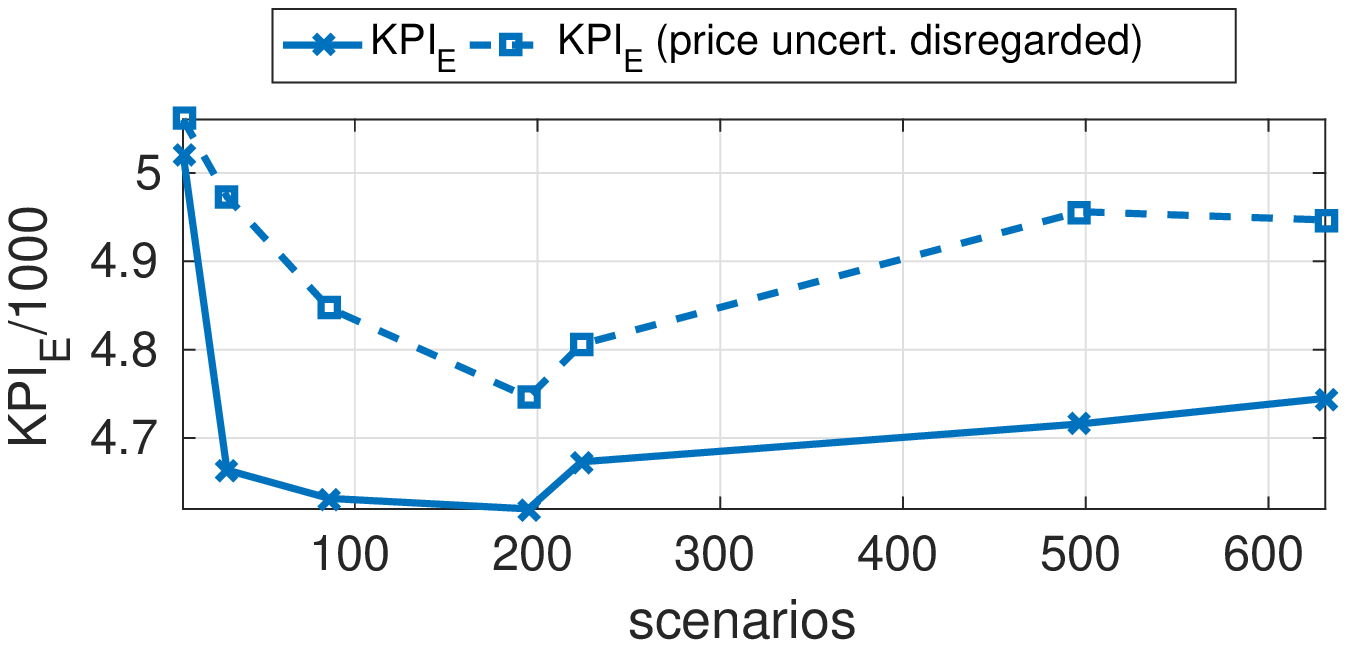}
 \includegraphics[width=0.98\columnwidth]{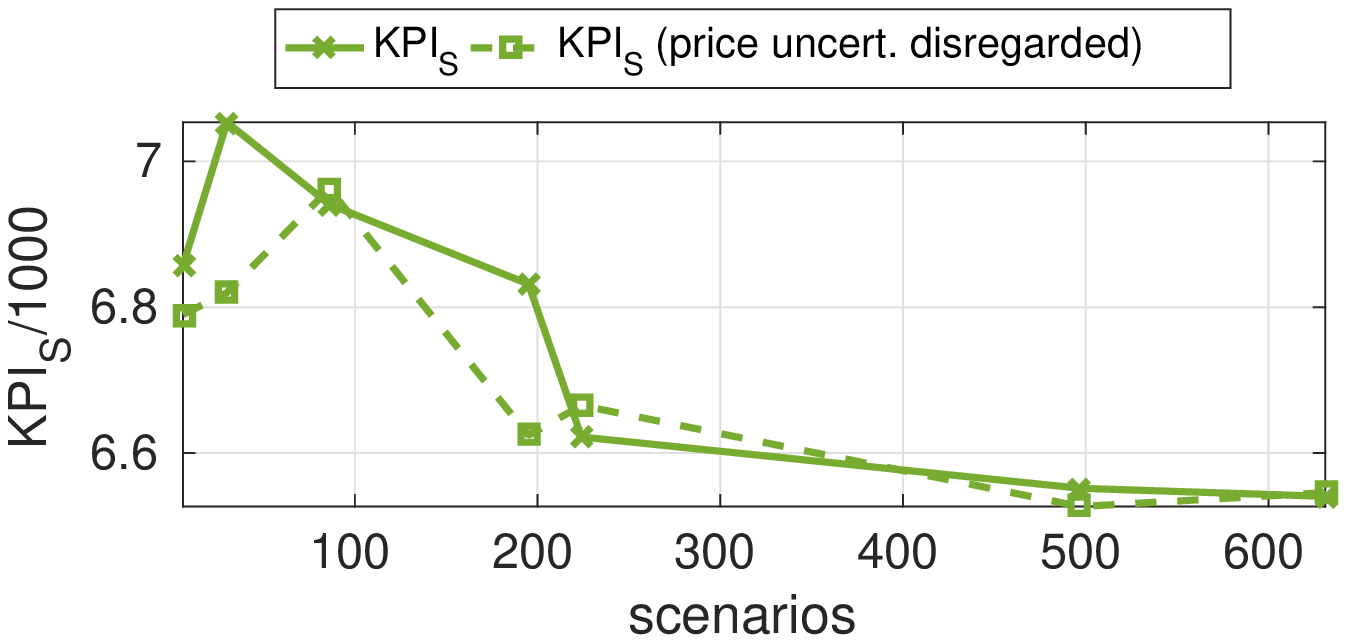}
 \caption{(Up) The economic index $\mathrm{KPI}_{E}$ {in euros (\euro)} as a function 
          of the number of scenarios considered in the SSMPC formulation,
          (Down) The safety index $\mathrm{KPI}_{S}$.
          The dashed lines correspond to SSMPC's which do not 
          consider the volatility in electricity prices.}
 \label{fig:KPIES}
\end{figure}

In order to justify the need for and underline the importance of 
taking into account the uncertainty associated with electricity 
prices, we performed two sets of simulations where in one 
case we disregarded that volatility (the SSMPC used only the 
nominal predictions of the electricity prices).
We may observe in Fig.~\ref{fig:KPIES} that the hourly operation 
of the network becomes more expensive by approximately $+5\%$,
therefore, the proposed stochastic control approach can lead to 
significant economic savings for the network operator.
{For example, for the case of $631$ scenarios, the operation 
cost when the price uncertainty is disregarded is \EUR{4947}/hour, whereas
with the proposed approach which takes into account the price uncertainty 
it drops to \EUR{4748}/hour --- a saving of \EUR{199}/hour, that is, a cost reduction of $4\%$.}

We may also see that as more scenarios are considered in the 
SSMPC formulation, the average cost of operation plummets at
around $194$ scenarios, where, however, the safety index is still 
high. In order to obtain a safer operation we need to pay the 
``cost of safe operation'': Indeed, at $631$ scenarios, the 
operation of the network is more expensive compared to the case of 
$194$ scenarios, but $\mathrm{KPI}_S$ is at a minimum.
This reflects a trade-off between the economic and safe operation 
of the network.

\begin{figure}[!]
 \centering
 \includegraphics[width=0.98\columnwidth]{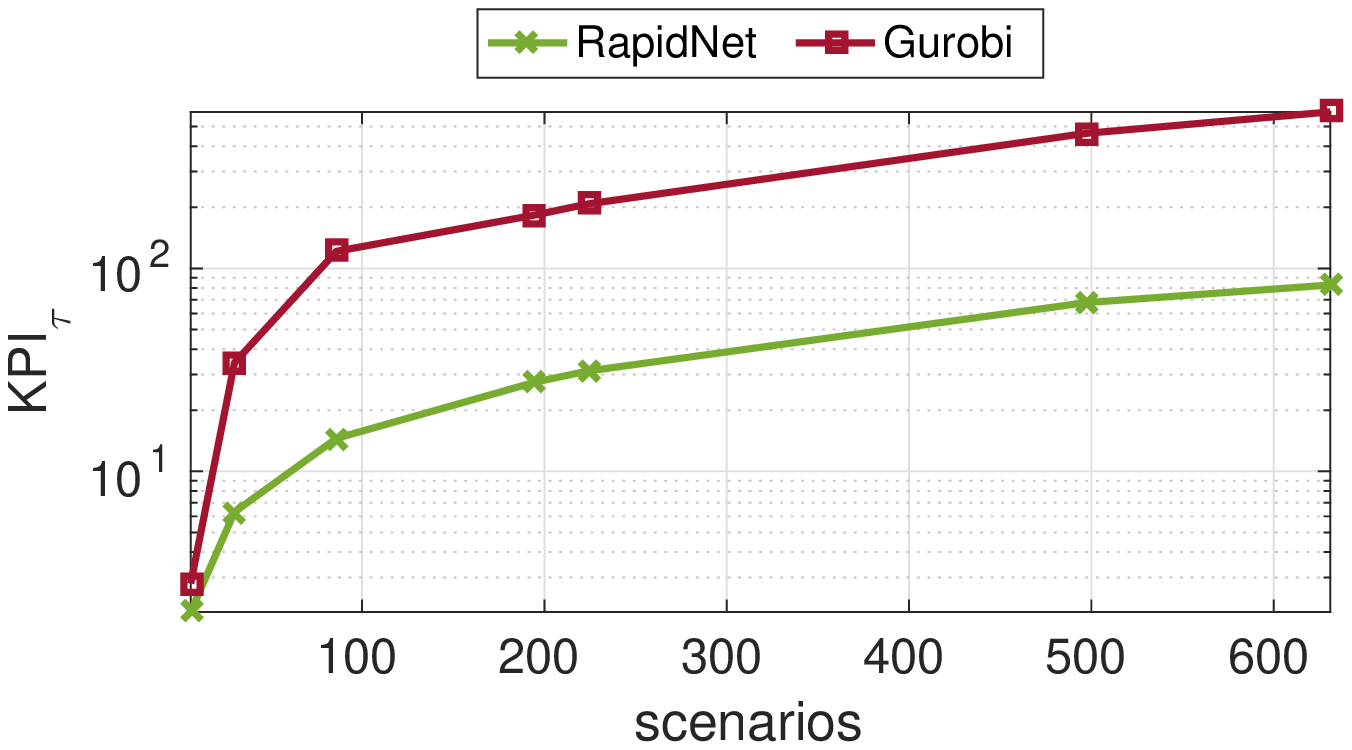}
 \caption{$\mathrm{KPI}_\tau$ for Gurobi (executed on an Intel Core i7-6600U 
  machine with  4 $\times$ 2.60GHz CPUs, 12GB RAM running 64-bit Ubuntu 16.04)
  and \texttt{\name} running on an NVIDIA Tesla C2075 GPU.}
 \label{fig:KPItau}
\end{figure}

At $631$ scenarios, the SSMPC optimization problem involves 
as many as $2,306,133$ primal variables and $3,126,960$
dual variables. \texttt{\name} solves it with $\mathrm{KPI}_{\tau} = \unit[82.7]{s}$
as compared to $\mathrm{KPI}_{\tau} = \unit[719]{s}$ for the popular commercial 
solver Gurobi v7.0 at the same level of accuracy (tolerance $5\cdot 10^{-2}$).
The solvers \textsc{cplex} and \textsc{mosek} were also tested, but Gurobi outperformed them
consistently.
As we may see in Fig.~\ref{fig:KPItau}, \texttt{\name} exhibits a lower complexity
index by approximately an order of magnitude.

\section{Conclusions and Future Work}
We presented the integrated software solution \texttt{\name} for the control of
drinking water networks which accounts for uncertainty both
in predicted water demand and in predicted electricity prices in the
day-ahead energy market. 

\texttt{\name} is a highly inter-operable software as it can be 
interfaced using JSON files which follow a standard API which is 
detailed in the software documentation. It can be combined with 
any forecaster as the controller does not need to know the mechanism
with which the forecasts are produced. {\texttt{\name} features 
a parser implemented in MATLAB which allows the conversion of EPANET 
\texttt{.inp} files to JSON files.} It is an open-source and free
software which is distributed under the terms of the GNU LGPL v3.0 licence
and can be downloaded at \url{https://github.com/GPUEngineering/RapidNet}.

In this paper we have shown via simulations that when a water network is operated
in the context of a volatile energy market, considerable savings can be 
obtained by using a predictor of the upcoming energy prices and taking 
into account the associated price-related uncertainty in the formulation 
of the scenario tree. Alongside, we advocate that \texttt{\name} can 
transform SSMPC from a powerful (but too complex) theoretical development
to control engineering practice and enabling the solution of very large 
SSMPC scale problems and their seamless integration into the control 
system of the water network.

This work focuses on flow-based water distribution networks such as
the network of the city of Barcelona. Pressure-based networks lead to 
optimization problems with nonconvex constraints. These can be approximated 
by solving a constraints satisfaction problem (CSP) as discussed in~\cite{ConPuiCem14}.
The problem can be then transformed into a convex SSMPC problem which 
can be solved by \texttt{\name}. {
The operator splitting concept has been proven to apply to sums of 
nonconvex functions which allows the solution of nonconvex optimal control
problems with input constraints~\cite{themelis2016forward}. This allows the 
application of recent developments in nonconvex optimization such 
as PANOC~\cite{panoc2017} for the solution of nonlinear MPC problems.}

% { According to~\cite{CST14}, ``heuristic control 
% laws with intermittent operator interventions are used today.''
% With this work we }

Future developments in \texttt{\name} will involve the implementation of
new faster parallelizable algorithms which make use of quasi-Newton 
directions based on our recent theoretical work~\cite{SamSopBemPat17a}
and the exploitation of multiple-GPU architectures. 
We shall introduce a library of available water network entities
from the literature as well as reported water demand and energy price models
for different water networks and energy markets.

%% The Appendices part is started with the command \appendix;
%% appendix sections are then done as normal sections
%% \appendix

%% \section{}
%% \label{}

%% If you have bibdatabase file and want bibtex to generate the
%% bibitems, please use
%%
%%  \bibliographystyle{elsarticle-num}
%%  \bibliography{<your bibdatabase>}

%% else use the following coding to input the bibitems directly in the
%% TeX file.
\section*{References}
{
\bibliographystyle{elsarticle-num}
\bibliography{biblio}
}

\end{document}